\documentclass[12pt, amssymb, amsmath]{article}

\begin{document}

\renewcommand{\thesection}{\arabic{section}}
\renewcommand{\thesubsection}{\thesection.\arabic{subsection}}

%-------------Lemmes, theoremes et autres---------------------------

%\newtheorem{lem}{Lemme}[section]
\newtheorem{lem}{Lemma}[section]
\newtheorem{propo}[lem]{Proposition}
\newtheorem{theo}[lem]{Theorem}
\newtheorem{rema}[lem]{Remark}
\newtheorem{remas}[lem]{Remarks}
\newtheorem{coro}[lem]{Corollary}
\newtheorem{defin}[lem]{Definition}
\newtheorem{hypo}[lem]{Hypoth\`ese}
\newtheorem{exem}[lem]{Exemple}
\newtheorem{conj}[lem]{Conjecture}

\newcommand{\real}{{\bf R}}
\newcommand{\reall}{{\bf R}^{p}}
\newcommand{\realll}{{\bf R}^{k}}
\newcommand{\mreal}{M\times{\bf R}^{p}}
\newcommand{\mreall}{M\times{\bf R}^{2k}}
\newcommand{\ri}{\rightarrow}
\newcommand{\tast}{T^{\ast}M}
\newcommand{\tastl}{T^{\ast}L}
\newcommand{\tastilde}{T^{\ast}\widetilde{M}}
\newcommand{\cd}{D_{\bullet}}
\newcommand{\cdk}{D_{\bullet-k}}
\newcommand{\cfx}{C_{\bullet}(f,\xi )}
\newcommand{\cald}{{\cal D}(M)}
\newcommand{\cm}{C_{\bullet}(M)}
\newcommand{\cmk}{C_{\bullet-k}(M)}
\newcommand{\cov}{\widetilde{M}}
\newcommand{\calm}{{\cal M}}
\newcommand{\cala}{{\cal A}}
\newcommand{\calc}{{\cal C}}
\newcommand{\calitate}{{\cal L}}
\newcommand{\caly}{{\cal Y}}
\newcommand{\gol}{\, \, \, \, \, \, \, \, }
\newcommand{\hatf}{\hat{f}}
\newcommand{\hatx}{\hat{x}}
\newcommand{\hatv}{\hat{v}}
\newcommand{\mn}{M^{n}}
\newcommand{\barm}{\bar{M}}
\newcommand{\bfs}{{\bf S}}
\newcommand{\piu}{\pi_{1}}
\newcommand{\pii}{\pi_{i}(M)}
\newcommand{\omi}{\Omega(L_{0},L_{1})}
\newcommand{\omibar}{\Omega(\bar{L}_{0},\bar{L}_{1})}
\newcommand{\omil}{\Omega(L,L)}
\newcommand{\omibarl}{\Omega(\bar{L},\bar{L})}

\begin{center}
{\large\bf  CONSTRAINTS ON 
EXACT LAGRANGIANS IN COTANGENT
BUNDLES OF MANIFOLDS FIBERED OVER THE CIRCLE}\\

\vspace{.3in}

\noindent Mihai DAMIAN\\
\noindent Universit\'e Louis Pasteur\\
IRMA, 7, rue Ren\'e Descartes,\\
67 084 STRASBOURG\\
e-mail : damian@math.u-strasbg.fr\\

\end{center}
\vspace{.4in}

\noindent{\bf Abstract\, :}\,  We give topological obstructions to the existence of a
closed exact Lagrangian submanifold $L\hookrightarrow \tast$, where $M$ is the total
space of a fibration over the circle. For instance, we show that $\piu(L)$ cannot be
the free product of two non-trivial groups and that the difference between the
number of generators and the number of relations in a finite presentation of
$\piu(L)$ is less than two. 
\vspace{.2in}

\section{Introduction}

Let $M^{n}$  be a closed connected  manifold  and $\tast$ its cotangent
bundle endowed with the standard symplectic structure $\omega_{M}\, =\,
d\lambda _{M}$, where $\lambda _{M}$ is the Liouville form
 $\lambda _{M}\, =\, \sum_{i} p_{i}\, dq_{i}$. Let $L^{n}\hookrightarrow \tast$ be an
exact Lagrangian submanifold, i.e. a submanifold such that $\lambda _{M}|{L}$ is an
exact $1$-form. 

The only known examples of exact Lagrangian submanifolds are the
graphs of functions $f:M\ri\real$ $$L_{f}:=\{(q,df_{q})\, |
\, q\, \in\, M\}$$
and their images by Hamiltonian vector flows.  The question of the
existence of other examples was first evoked by V.I. Arnold in his survey "First steps
in symplectic topology" \cite{Arn}. It is far from being solved. A positive answer
was given by R. Hind in the case $L=M={\bf S}^{2 }$. The other related
results which were proved up to now are  
topological obstructions to the existence of exact Lagrangian embeddings
$L\hookrightarrow\tast$. We summarize them in the statement below :\\

\noindent{\bf Theorem 0} \,  {\it Let $M$ be a closed manifold and $L\hookrightarrow \tast$
an exact Lagrangian embedding of a closed manifold $L$. Denote by $p$ the projection
of $L$ on the base space $M$. Then we have : \\
a) If $L$ and $M$ are orientable, then $\chi(L)\, =\, deg^{2}(p)\chi(M)$. If $L$ and
$M$ are not orientable the same equality is valid modulo $4$.\\
b)  The index $[\pi_{1}(M)\, :\, p_{\ast}(\pi_{1}(L))]$ is finite. \\
c) If $M$ is simply connected then $L$ can not be aspherical (i.e.
Eilenberg-Mac Lane). \\
d) If $M$  and $L$  are spin and $L$ has a vanishing Maslov class, 
then $H^{\ast}(L,K)\approx 
H^{\ast}(M,K)$, for any field $K$  with $char(K)\neq 2$.\\
e) If $M$ is simply connected then $p^{\ast}:H^{2}(M)\ri H^{2}(L)$ is injective and 
the index $[\pi_{2}(M)\, :\, p_{\ast}(\pi_{2}(L))]$ is finite. } 

\vspace{.2in}

The statement 0.a was proved by M. Audin in \cite{Au}, 0.b was
proved by F. Lalonde and J-C Sikorav in \cite{LS} and  0.c is a result of C.
Viterbo \cite{Vi1} (see also \cite{Vi}). More recently, 0.d was proved 
independently by 
K. Fukaya, P. Seidel and I. Smith \cite{FSS1}, \cite{FSS2} and D. Nadler \cite{N}. 
For $M={\bf S}^{n}$ and $L$ simply connected this was proved previously by
P. Seidel \cite{Sei} 
and by L. Buhovsky \cite{Buh}. 

The assertion 0.e was proved by A. Ritter \cite{Rit}. 
The author uses some techniques coming from the Novikov homoloy theory which we will
also do below.   

  The aim of this paper is to provide other obstructions in the case where $M$ is a
total space of a fibration over the circle. Let us state our main results~: 

\begin{theo} \label{main}  Let $M^{n\geq 3}$ be a closed manifold 
which is the total space of a
fibration over ${\bf S}^{1}$ and let $L\hookrightarrow \tast$ be an exact Lagrangian
embedding of a closed manifold $L$. Then we have :\\
a) Let $<g_{1}, g_{2}, \ldots, g_{p}\, |\, r_{1}, r_{2}, \ldots, r_{q}>$ be an
arbitrary presentation of the fundamental group $\piu(L)$. Then $p-q\leq 1$. \\
b) The fundamental group  
$\pi_{1}(L)$ is not isomorphic to the free product $G_{1}\ast G_{2}$ of two
non-trivial (finitely presented) groups.
\end{theo}

Here are some exemples of non-embedding statements which can be inferred from
our result :

\begin{coro}\label{coro} Let $P,Q, L$ be closed manifolds and suppose 
 that $\piu(P)$ is finite. \\
a) Suppose that 
$\chi(L)\neq 0$. Then there is no exact Lagrangian
embedding $L\times P\hookrightarrow 
T^{\ast}(Q\times{\bf S}^{1})$. 

In particular, let  $\Sigma_{g}$ be a (non necessary orientable) surface of genus
 $g\geq 2$. Then there is no exact Lagrangian embedding of $\Sigma_{g}\times P$ into
$T^{\ast}(Q\times{\bf S}^{1})$. More generally, for surfaces $\Sigma_{{g}_{i}}$ as
above there is no exact Lagrangian embedding
$$\Sigma_{g_{1}}\times\Sigma_{g_{2}}\times\cdots\times\Sigma_{g_{k}}\times 
P\hookrightarrow 
T^{\ast}(Q\times{\bf S}^{1}).$$
b)
  Let $L^{n\geq 4}$ be the connected sum $L_{1}\# L_{2}$ of two closed manifolds.
 Then there is no exact Lagrangian
embedding $L\times P \hookrightarrow 
T^{\ast}(Q\times{\bf S}^{1})$ unless one of the $L_{i}$ is a simply connected ${\bf
Z}/2$-homology sphere.\\
c) Suppose that there is an exact Lagrangian embedding $$L\times
T^{l}\hookrightarrow T^{\ast}(T^{m}\times Q),$$ where $T^{k}$ is the
$k$-dimensional torus and  $m>l$. Then $L$ satisfies the conditions a, b of
\ref{main}.            \end{coro}

\subsection{Idea of the proof} 

Let $f:M\ri\bfs^{1}$ be a fibration. The closed $1$-form $\alpha=f^{\ast}d\theta$ has
no zeroes. Let $L$ be an exact embedding into $\tast$. 
Consider the Lagrangian isotopy $$L_{t}\, =\, L+t\alpha.$$
It follows that $L_{t}\cap L=\emptyset$ for $t$ large enough. The Lagrangian manifolds
$L_{t}$ are not exact but they satisfy 
$$\omega_{M}|_{\pi_{2}(\tast,L_{t})}=0$$ just like an exact Lagrangian manifold.
 Under this
hypothesis one can define a Floer-type complex  $C_{\bullet}(L,L_{t})$, which is spanned
by the intersection points $L\cap L_{t}$. Therefore, this complex vanishes for $t>>0$. 

It turns out that the homology of this
complex is isomorphic to the Novikov homology $H_{\ast}(L,p^{\ast}u)$, where
$u\in H^{1}(M,{\bf Z})$ is the cohomology class of $\alpha$ and $p:L\ri M$ is the
projection. In particular it is independent of $t$. It follows : 

\begin{theo}\label{novi} $H_{\ast}(L,p^{\ast}u)=0$\end{theo}

In order to prove \ref{main} one has to argue in the following way
 : suppose that \ref{main}.a
is false. Then one can show  that the Novikov homology $H_{\ast}(L,v)$ does not vanish
for any $v\in H^{1}(L,\real)$, contradicting thus \ref{novi}. A similar argument works
for the proof of \ref{main}.b. 

The paper is organized as follows. In Section 2 we recall the definition and the main
properties of Novikov homology. We prove the non-vanishing results needed in the above
argument.  In Section 3 we define the Floer complex $C_{\bullet}(L, L_{t})$. Finally, in
Section 4 we establish \ref{novi} and prove \ref{main} and \ref{coro}. 

\section{Novikov theory}
\subsection{Definition of Novikov homology} 

Let $u\in H^{1}(L;\real)$.
Denote by $\Lambda$ the ring   $ {\bf Z}/2\, [\pi _{1}(L)]$ and by $ \hat{\Lambda}$ the
ring of formal series ${\bf Z}/2\, [[\pi_{1}(L)]]$.
Consider a $CW$-decomposition of $L$
which we lift  to the universal
cover $\widetilde{L}$. We get a $\Lambda$-free complex $C_{\bullet}(\widetilde{L})$
  spanned by
(fixed lifts of) the cells of the triangulation of $L$.  

We define now {\it the completed ring} $\Lambda_{u}$  :
$$\Lambda_{u}:=\, \left\{\lambda\, =\, \sum n_{i} g_{i}\, \in\,
\hat{\Lambda}\, \, |\, \, g_{i}\in \pi_{1}(L), \, \, n_{i}\in {\bf Z}/2, \, \, \,
u(g_{i})\, \rightarrow\,  +\infty\right\}.$$
The convergence to $+\, \infty$ means
 here that for all $A\, >\, 0$,  $u(g_{i})\, <\, A$ only for a finite
number of $g_{i}$ which appear with a non-zero coefficient in  the sum  $\lambda$.\\

\begin{rema}\label{inverse}  Let $\lambda\, =\, 1+\sum n_{i} g_{i}\, \in \, \Lambda_{u}$ 
where $u(g_{i})>0$ for all $i$.
Then $\lambda$ is invertible in $\Lambda_{u}$. Indeed, if we 
denote by $ \lambda_{0}\, =\,
\sum n_{i} g_{i}$ then it is easy to check that $\sum_{k\geq 0}(-\lambda_{0})^{k}$ 
is an element of
$\Lambda_{u}$ and it is obvious that it is the inverse of $\lambda$.\end{rema}

\noindent{\bf Definition} Let $C_{\bullet}(L,u)$ be the  
$\Lambda_{u}$-free complex
$\Lambda_{u}\otimes_{\Lambda}C_{\bullet}(\widetilde{L})$. 
The Novikov homology $H_{\ast}(L,u)$ is
the homology of the complex $C_{\bullet}(L,u)$.\\
\\
\begin{rema}\label{coefficients}  We may define in a similar way the Novikov 
homology with $\bf Z$
coefficients. As we want to compare it to Floer homology and the latter is
defined for ${\bf Z}/2$ coefficients we used ${\bf Z}/2$ in the definitions above.
\end{rema}

  Now we prove :

\begin{propo}\label{nonzero} Let $L$ and $u$ be as above. \\
a)  Let 
 $<g_{1}, g_{2}, \ldots, g_{p}\, |\, r_{1}, r_{2}, \ldots, r_{q}>$ be a
 presentation of the fundamental group $\piu(L)$ which satisfies $p-q\geq 2$. Then, if
$u\neq 0$ we have   
$H_{1}(L,u)\neq 0$. \\
b) Suppose that $\piu(L)=G_{1}\ast G_{2}$, none of the $G_{i}$ being trivial. Then, if
$u\neq 0$ we have   
$H_{1}(L,u)\neq 0$.
\end{propo}

\noindent\underline{Proof}

The presentation of $\piu(L)$ yields a $CW$ decomposition of $L$ with one single
zero-cell, $p$ one-cells and $q$ two-cells. Lifting it to the universal cover we see
that the complex $C_{\bullet}(\widetilde{L})$ ends like follows :
$$  \cdots \ri\Lambda^{q}\stackrel{\delta_{2}}{\ri}
\Lambda^{p}\stackrel{\delta_{1}}{\ri}\Lambda\ri 0.$$
If $\{e_{i}\}_{i=1,\ldots,p}$ is the basis of $C_{1}=\Lambda^{p}$ given by the
$1$-cells
 and $\{e\}$ is the basis of $C_{0}=\Lambda$ given by the single $0$-cell, it is easy
to see that the differential $\delta_{1}$ satisfies $\delta_{1}(e_{i})=(1-g_{i})e$.\\
\\
a) Since $u\neq 0$, $u(g_{i})\neq 0$ for some $i$, we have  that
$1-g_{i}$ is invertible in $\Lambda_{u}$ (its inverse is $1+g_{i}+g_{i}^{2}+\cdots$ if
$u(g_{i})>0$, resp. $-g_{i}(1+g_{i}^{-1}+g_{i}^{-2}+\cdots)$ if $u(g_{i})<0$).
In particular in the tensored complex $$(1)\gol \cdots
\ri\Lambda_{u}^{q}\stackrel{\delta_{2}}{\ri}
\Lambda_{u}^{p}\stackrel{\delta_{1}}{\ri}\Lambda_{u}\ri 0.$$
 the dimension of $Ker(\delta_{1})$ is $p-1$. Since $p-1>q$, it follows
that $Im(\delta_{2})\neq Ker(\delta_{1})$, so $H_{1}(L,u)\neq 0$ as claimed. \\
\\b) Let $\{f_{j}\}_{j=1\ldots, q}$ the basis of $\Lambda^{q}$ defined by the two-cells
corresponding to the relations $r_{1}, r_{2},\ldots, r_{q}$. The matrix of
$\delta_{2}$ with respect to $\{f_{j}\}$ and $\{e_{i}\}$ is given by the Fox
derivatives $\partial r_{j}/\partial g_{i}$ \cite{Fox}. These derivatives are defined
by the following formulas :
$$\frac{\partial g_{i}}{\partial g_{i}}=1, \gol, \frac{\partial g_{i}^{-1}}{\partial
g_{i}}=-g_{i}^{-1}, \gol \forall i=1,\ldots p,$$
$$\frac{\partial(rr')}{\partial g_{i}}=\frac{\partial r}{\partial
g_{i}}+r\frac{\partial r'}{\partial g_{i}}, \gol, \forall i=1,\ldots p,$$
where $r$ and $r'$ are words written with the letters $g_{i}^{\pm 1}$. 

Now suppose that $\piu(L)=G_{1}\ast G_{2}$ and consider (for $k=1,2$) 
finite presentations of
$G_{k}$ with $p_{k}$ generators and $q_{k}$ relators. Denote  by $\delta_{1}^{k}$ and
$\delta_{2}^{k}$ the differentials of the complex $(1)$ corresponding to these finite
presentations. Then, for some $u:\piu(L)\ri\real$, using the
 definition of the maps $\delta_{1}$ and $\delta_{2}$, we find that the
complex $(1)$ writes : 
$$\cdots
\ri\Lambda_{u}^{q_{1}}\oplus\Lambda_{u}^{q_{2}}
\stackrel{\left(\begin{array}{cc} \delta_{2}^{1} & 0\\
0& \delta_{2}^{2}\end{array}\right)}{\longrightarrow}
\Lambda_{u}^{p_{1}}\oplus \Lambda_{u}^{p_{2}}
\stackrel{(\begin{array}{cc} \delta_{1}^{1}& \delta_{1}^{2}\end{array})
 }{\longrightarrow}\Lambda_{u}\ri 0.$$

Suppose now that $u\neq 0$ and, without restricting the generality, that
$u|_{G_{1}}\neq 0$. As above the map $\delta_{1}^{1}$ is then surjective. This
implies that for any $a\in \Lambda_{u}^{p_{2}}$ there is an element $b\in
 \Lambda_{u}^{p_{1}}$ such that $(b,a)$ belongs to the kernel of  
$(\begin{array}{cc} \delta_{1}^{1}& \delta_{1}^{2}\end{array})$.

If $H_{1}(L,u)=0$, it follows that
$\delta_{2}^{2}$ is an epimorphism. Now the sequence 
$$\Lambda_{u}^{q_{2}}\stackrel{\delta_{2}^{2}}{\ri}
\Lambda_{u}^{p_{2}}\stackrel{\delta_{1}^{2}}{\ri}\Lambda_{u}\ri 0$$
is exact and therefore $\delta_{1}^{2}=0$. But this is impossible unless $G_{2}=0$ and
the proof is finished. \hfill $\diamond$ 

\subsection{Morse-Novikov theory}

We recall in this subsection the relation between Novikov homology and 
closed $1$-forms. 
Let $\alpha$ be a closed generic $1$-form in the class $u\in H^{1}(L,\real)$.
"Generic" means here that the zeroes of $\alpha$ are of Morse type.
  Let $\xi$ be the
gradient of $\alpha$ with respect to some generic metric on $L$. For every zero
 $c$ of $\alpha$ we fix a lift $\tilde{c}$ of $c$ in the universal cover
$\widetilde{L}$.  We can define then a
complex $C_{\bullet}(\alpha,\xi)$ spanned by the zeroes of $\alpha$ and graded by the
Morse index :
the "incidence number" $[d,c]$ for two zeroes
of consecutive Morse indices is the (possibly infinite) sum
 $\sum n_{i} g_{i}$ where $n_{i}$ is the algebraic number of
flow lines which join $c$ and $d$ and which are covered by a path in $\widetilde{L}$
 joining
$g_{i}\tilde{c}$ and $\tilde{d}$. It turns out that this incidence number belongs to
$\Lambda_{u}$, so $C_{\bullet}(\alpha,\xi)$ is actually a $\Lambda_{u}$-free complex.

 The fundamental property of the Morse-Novikov  theory was proved by S.P. Novikov in
\cite{Nov} and generalized by J.-C. Sikorav in \cite{Sikt}. The statement is :
\begin{theo}\label{novikov}
For any generic couple $(\alpha,\xi)$ as above, the
homology of the complex $C_{\bullet}(\alpha,\xi)$ is isomorphic to
$H_{\ast}(L,u)$.\end{theo}

\begin{rema}\label{novikovgeneral} We may define in a similar way the Novikov homology 
$H_{\ast}(L,u)$ associated to a
covering $\pi:\bar{L}\ri L$ which satisfies $\pi^{\ast}(u)=0$ and we also may define
$C_{\bullet}(\alpha,\xi)$ using the covering $\bar{L}$ instead of $\widetilde{L}$.
The statement \ref{novikov} holds in this more general setting. 
 \end{rema}

  An easy consequence of \ref{novikov} is the following statement : 

\begin{propo}\label{Kunneth}
 Let $L_{1}$, $L_{2}$ be closed manifolds and let
$u\in H^{1}(L_{1},\real)\subset H^{1}(L_{1}\times L_{2},\real)$.
 Consider the Novikov homology $H_{\ast}(L_{1}\times L_{2},
u)$ associated to $u$ and to a covering $(\pi, Id):\bar{L}_{1}\times L_{2}\ri
L_{1}\times L_{2}$. 
 Then $$H_{\ast}(L_{1}\times L_{2},
u)\, \approx\,   H_{\ast}(L_{1},u)\otimes_{{\bf Z}/2} H_{\ast}(L_{2};{\bf
Z}/2).$$
In particular, 
 $$H_{\ast}(L_{1}\times L_{2},
u)=0\, \Longrightarrow  H_{\ast}(L_{1},u)=0.$$
\end{propo}

\noindent\underline{Proof}

Note first that the homologies $H_{\ast}(L_{1}\times L_{2},
u)$ and  $H_{\ast}(L_{1},u)$ are defined using the same Novikov ring  $\Lambda_{u}$.
Take a generic pair $(\alpha_{1},\xi_{1})$ associated to $u$ on $L_{1}$ and a
generic pair $(df_{2},\xi_{2})$ on $L_{2}$. One can easily see that the complex
$C_{\bullet}(\alpha_{1}+df_{2}, \xi_{1}+\xi_{2})$ is isomorphic as a 
$\Lambda_{u}$-complex to the tensor product
$C_{\bullet}(\alpha_{1},\xi_{1})\otimes_{{\bf Z}/2}C_{\bullet}(df_{2},\xi_{2})$. 
 By comparing their homologies using the Kunneth formula we get
$$H_{\ast}(L_{1}\times L_{2},
u)\, \approx\,   H_{\ast}(L_{1},u)\otimes_{{\bf Z}/2} H_{\ast}(L_{2};{\bf
Z}/2)$$
as claimed \hfill
$\diamond$ 

\vspace{.2in} 

We end this section by the following trivial remark :

\begin{rema}\label{nonze} Let $L$ be a manifold with $\chi(L)\neq 0$. Then the
Novikov homology, defined for any  covering $\pi:\bar{L}\ri L$, as above, does not
vanish.\end{rema}

 Indeed, the complex $C_{\bullet}(L,u)$ has the same Euler characteristic as $L$ and if the 
 Novikov homology $H_{\ast}(L,u)$ vanishes, this Euler characteristic must be zero.

 \section{The Floer complex}

In \cite{Leo}  Hong Van Le and Kaoru
 Ono defined a Floer complex spanned by the $1$-periodic
orbits of a symplectic {\it non Hamiltonian} vector field $X_{t}$ on a symplectic
manifold $W$. They showed that, if $W$ is monotone, then this homology equals the
Novikov homology associated to the cohomology class $Cal(\phi_{t})$ ; by
$Cal(\phi_{t})$ we denote the Calabi invariant associated to the symplectic flow
$(\phi_{t})$ of $X_{t}$, defined by the formula 
$$Cal(\phi_{t})\, =\, [\int_{0}^{1}\beta_{t}dt],$$ where $\beta_{t}$ is the
symplectic dual of $X_{t}$ for $t\in [0,1]$. This is actually the image of the isotopy
$(\phi_{t})$ by the Flux morphism. By definition, the integral of $Flux(\phi_{t})$ over
a loop $c:{\bf S}^{1}\ri W$ is the integral of the symplectic form over the cylinder
$\phi_{t}(c)$.

 Later, M. Pozniak \cite{Poz}
 computed the  (Lagrangian) Floer homology
$HF(L_{1},L_{2})$ in the case where $L_{1}$ is the zero-section of $\tast$ and 
   $L_{2}=\phi(L_{1})$, $\phi$ being the time one map of a symplectic 
(non-Hamiltonian) isotopy. 

 Consider a closed exact Lagrangian manifold 
$L\hookrightarrow \tast$ and denote by $L_{t}$ the image of $L$ through a
symplectic  isotopy $(\phi_{t})$ on $\tast$. 
Denote by $u\in H^{1}(M;\real)$ the class $$
[\phi^{\ast}_{1}\lambda_{M}-\lambda_{M}]\, \in\, H^{1}(\tast;\real)\approx
H^{1}(M;\real).$$ It is not difficult to prove that this
 is actually the Calabi invariant  $Cal(\phi_{t})$ for 
$W=\tast$. 
  The goal
of this section is the following 

\begin{theo}\label{conclusion}
Suppose that $L$ and $\phi_{1}(L)$ are transverse and that $u$ is rational (this
means that the image of the morphism $u:\piu(M)\ri \real$ is cyclic). 
Denote by $p^{\ast}u$ the
composition $u\circ p :\piu(L)\ri{\bf Z}$. 
There is a free $\Lambda_{p^{\ast}u}$-complex 
$C_{\bullet}(L,\phi_{t})$ spanned by the
intersection points $L\cap\phi_{1}(L)$ and whose homology only depends on $L$ and on
$u$. 
\end{theo}

\begin{rema}\label{Tischler}
The restriction $Im(u)\approx {\bf Z}$ is not strong. Indeed, the isotopy 
$(\phi_{t})$
is the flow of a vector field $X_{t}$ whose symplectic dual is a family of 
 closed one forms $\beta_{t}$. If $\beta$ belongs to the class $u=Cal(\phi_{t})$,
 a well known result of D. Tischler \cite{Ti} asserts that there exists a closed
one-form $\beta'$, arbitrarily closed to $\beta$ and whose cohomology class 
$[\beta']$
is rational. This implies that  given a simplectic isotopy $(\phi_{t})$, there is an
 symplectic isotopy $(\phi'_{t})$ arbitrarily closed to it and with rational Calabi
invariant : one can take $(\phi'_{t})$ defined by the family of closed one forms
$\beta_{t}+\beta'-\beta$. 
\end{rema}

\noindent\underline{Plan of the proof of \ref{conclusion}}

The complex $C_{\bullet}(L,\phi_{t})$ will be the result of a version of Novikov homology theory on the
(infinite dimensional) space of paths joining $L$ and $\phi_{1}(L)$.  After some preliminary results
proved in the next subsection we define in Subection \S3.2 a action one-form $\nu$ whose 
zeroes are in one-to-one
correspondence with the intersection points of $L$ and $\phi_{1}(L)$. In the following 
subsection we consider
the gradient of $\nu$ with respect to a metric defined by a family of almost complex structures on
$\tast$. The "flow lines" of this gradient vector fields will be
actually holomorphic strips with boundary on $L\cup\phi_{1}(L)$. 
The results of \S3.4 enable us to use this Floer-Novikov setting to define a free complex over
the Novikov ring $\Lambda_{u}$ which is spanned by the points of $L\cap\phi_{1}(L)$ (Subsection \S3.5). 
 In \S3.6 we prove
that the homology of this
complex only depends on $L$ and $u$. Finally, we show in \S3.7 how the whole construction  can be adapted
in order to get a free complex over the Novikov ring $\Lambda_{p^{\ast}u}$, spanned by $L\cap\phi_{1}(L)$
and  whose homology again only depends on $L$ and on $u$.

Let us now explain in detail how one defines
the Floer complex $C_{\bullet}(L,\phi_{t})$. 

\subsection{ Preliminary results}

Let $(\phi_{t})$ a symplectic isotopy as above and 
denote by $u\in H^{1}(M;\real)$ the class $Cal(\phi_{t})$
We prove the following lemma~:

\begin{lem}\label{Calabi} There is a symplectic isotopy $(\psi_{t})$ on $\tast$ such
that $\psi_{1}|_{L}=\phi_{1}|_{L}$
 which is defined  by (the symplectic dual of) $\beta=\alpha+d H_{t}$, where $\alpha\in u$ is
a closed $1$-form on $M$ and $H:\tast\times[0,1]\ri\real$ has compact support.
\end{lem}

\noindent\underline{Proof}

The $1$-form $\phi_{t}^{\ast}\lambda_{M}-\lambda_{M}$ is
closed on $\tast$. So, we may write $$ \phi_{t}^{\ast}\lambda_{M}-\lambda_{M}\, =\,
\alpha_{t}+dG_{t},$$ where $G:\tast\times [0,1]\ri\real$ and the closed $1$-form
$\alpha_{t}$ is defined on $M$. 

 Define a symplectic isotopy $\Gamma_{t}:\tast\ri\tast$ by $$\Gamma_{t}(p,q)\, =\,
(p-\alpha_{t}(q),q).$$
Then, obviously, $\Gamma_{t}^{\ast}\lambda_{M}-\lambda_{M}=-\alpha_{t}$ ; Therefore,
for $\chi_{t}=\Gamma_{t}\circ\phi_{t}$ 
$$\chi_{t}^{\ast}\lambda_{M}=\phi_{t}^{\ast}(\lambda_{M}-\alpha_{t})=
\lambda_{M}+\alpha_{t}+dG_{t}-\phi_{t}^{\ast}\alpha_{t}=\lambda_{M}+dK_{t}$$
for a smooth $K:\tast\times[0,1]\ri\real$, since $\phi_{t}^{\ast}\alpha_{t}$ and 
$\alpha_{t}$ are cohomologous. 

 Using the Lie derivative, one obtains easily that $(\chi_{t})$ is a Hamiltonian
isotopy. We want the isotopy $(\chi_{t})$ 
 (and in particular the function $K$) to be compactly supported. Since $L$ is compact,
we may suppose that it is true
 and keep the relation $\chi_{t}=\Gamma_{t}\circ\phi_{t}$ valid on
$L$. In other words, we have $$\chi_{t}=\Gamma_{t}\circ\widetilde{\phi}_{t}, $$
where $\widetilde{\phi}_{t}|_{L}=\phi_{t}|_{L}$ and $\chi_{t}$ is compactly supported. 

  On the other hand $\Gamma^{-1}_{1}(p,q)=(p+\alpha_{1}(q),q)$ is the time one of the
symplectic isotopy $\widetilde{\Gamma}_{t}(p,q)=(p+t\alpha_{1}(q),q)$, so
$\widetilde{\phi}_{1}$ is the time one of $$\psi_{t}=\widetilde{\Gamma}_{t}\chi_{t}.$$

Outside a compact set $\psi_{t}(p,q)=\widetilde{\Gamma}_{t}(p,q)=(p+t\alpha_{1}(q),q)$
 is defined by $\alpha_{1}$. It
follows that $(\psi_{t})$ is defined by $\alpha_{1}+dH_{t}$ for some smooth, compactly
supported Hamiltonian 
$H:\tast\times[0,1]\ri\real$, and $\psi_{1}$
coincides with $ \phi_{1}$, when restricted to $L$, as claimed.
 Moreover, $[\alpha_{1}]$
is obviously the class $Cal(\phi_{t})$ in $H^{1}(M;\real)$. The proof is finished. 

\hfill $\diamond$

From now on, all the symplectic isotopies we consider have the property of Lemma
\ref{Calabi} and rational Calabi invariant. In the following we consider the lift of such an isotopy 
to the universal
cover of $\tast$.

  \begin{rema}\label{epi} Let $j:L\hookrightarrow\tast$ be a  
closed exact Lagrangian submanifold
and let $p:L\ri M$ be the
projection on the base space.  Without restricting the generality, we may add the
hypothesis $\piu(\tast,L)=0$ in the statement of \ref{main}. Indeed, by
Theorem 0.b, we know that the index $[\piu(M):p(\piu(L))]$ is finite. Consider the
finite cover $\bar{M}\ri M$ which corresponds to  the subgroup $p(\piu(L))$. Then, it
is easy to show that there is an exact Lagrangian embedding $L\hookrightarrow
T^{\ast}\bar{M}$ which is a lifting of $L\hookrightarrow\tast$ (the definition of this
lifting is similar to the one in 2.2 of
\cite{LS}). This embedding induces an epimorphism
$\pi_{1}(L)\ri\pi_{1}(\bar{M})$. Since $\bar{M}$ is still a total space of a fibration
over the circle, we may prove \ref{main} for $\bar{M}$ instead of $M$ in order to get
the desired obstructions on $L$. 
\end{rema}

Consider now the universal cover $\pi:\widetilde{M}\ri M$ and the induced projection
$\tilde{\pi}:T^{\ast}\widetilde{M}\ri\tast$.  
Denote by $y\mapsto g y$ the diffeomorphism of $\tastilde$ defined by the 
action of $g$ on $\piu(\tast)\approx\piu(M)$. Since this is a right action, one should keep in
mind that $g'(g''y)=(g''g')y$.  Let $K$ be the kernel of the epimorphism
$p:\piu(L)\ri\piu(M)$ and let $\bar{\pi}:\bar{L}\ri L$
 be the cover of $L$ associated to $K$. According to \ref{epi} this is a $\piu(M)$-covering, so
$\bar{L}$ is not compact.  We
prove 
\begin{lem}\label{cover}
There exists an exact Lagrangian embedding 
$\Psi: \bar{L}\hookrightarrow \tastilde$ which is a
lifting of $L\hookrightarrow \tast$. Moreover,
 we have  $\Psi(\bar{L})=\tilde{\pi}^{-1}(j(L))$
 and  
 for any $g\in\piu(M)$, 
$\bar{x}\in\bar{L}$,  $$\Psi(g\bar{x})=g\cdot
\Psi(\bar{x}).$$
\end{lem}

\noindent\underline{Proof}
 
 Consider the pullback of $\tilde{\pi}:\tastilde\ri \tast$ by the 
embedding $j:L\ri \tast$ i.e. the restriction to $L$ of the
covering $\tastilde\ri\tast$. Since $p:\piu(L)\ri\piu(M)$ is an epimorphism, this space
is connected.  It follows that  this  covering of $L$  is
isomorphic to $\bar{L} \ri L$. We keep the same notation 
$\bar{L}$ and consider the canonical
maps $\Psi:\bar{L}\ri\tastilde$ and $\bar{\pi}:\bar{L}\ri L$. 
 Note that, for $g\in \piu(M)$ we have
$$\Psi(g\bar{x})=g\cdot\Psi(\bar{x}),$$ as claimed. 

Using the commutative diagram 
$$\begin{array}{lcr}
\bar{L}&\, \, \stackrel{\Psi}{\ri}\, \, &\tastilde\\
\\
\downarrow\, \bar{\pi}&\, &\downarrow\, \tilde{\pi}\\
\\
L&\, \, \stackrel{j}{\ri}\, \, &\tast
\end{array}
$$ 
one easily checks that $\Psi$ is an exact Lagrangian embedding. The equality
$\Psi(\bar{L})=\tilde{\pi}^{-1}(j(L))$ is an obvious consequence of the definition of
the pullback. 

\hfill $\diamond$

\vspace{.2in} 

\noindent{\bf Notations} 
We use the same notation $\bar{L}\subset\tastilde$  for the Lagrangian submanifold which is the
image 
$\Psi(\bar{L})$ of the embedding 
constructed above. To keep the notations uniform, we will denote by $\bar{a}$ a lift to $\tastilde$
of an object $a\in \tast$. \\

 Consider now a symplectic isotopy $(\phi_{t})$ on $\tast$
 which is defined by $\alpha +dH_{t}$,
as in \ref{Calabi}. The following result is straightforward : 
\begin{lem}\label{fibres} The symplectic isotopy $(\phi_{t})$ lifts to a Hamiltonian
isotopy $(\bar{\phi}_{t})$ on $\tastilde$. 

  Moreover, if we denote $L_{t}=\phi_{t}(L)$ and
 $\bar{L}_{t}=\bar{\phi}_{t}(\bar{L})$, then 
$$\bar{L}\cap\bar{L}_{1}\, =\, \bigcup_{x\in
L\cap L_{1}}\tilde{\pi}^{-1}(x).$$ \end{lem}

\noindent\underline{Proof} 

Let $\alpha+dH_{t}$ be the closed $1$-form on $\tast$ whose symplectic dual $X_{t}$
 defines $(\phi_{t})$. Take
its pullback $\tilde{\pi}^{\ast}\alpha + d(H_{t}\circ\tilde{\pi})$ on $\tastilde$. The flow of
its symplectic dual $\bar{X}_{t}$ 
defines a symplectic isotopy $\widetilde{\phi}_{t}$ which is
actually Hamiltonian since $\tastilde$ is simply connected. It is obvious that
$\tilde{\pi}_{\ast}(\bar{X}_{t})=X_{t}$, which immediately implies that
$(\bar{\phi}_{t})$ is a lift of $(\phi_{t})$. We have thus a commutative diagram 
$$\begin{array}{lcr}
\tastilde&\, \, \stackrel{\bar{\phi }_{t}}{\ri}\, \, &\tastilde\\
\\
\downarrow\, \tilde{\pi}&\, &\downarrow\, \tilde{\pi}\\
\\
\tast&\, \, \stackrel{\phi_{t}}{\ri}\, \, &\tast
\end{array}
$$ 

 Now, following \ref{cover}, we have $\bar{L}=\tilde{\pi}^{-1}(L)$, and using the
above diagram $$\bar{L}_{t}=\bar{\phi}_{t}(\bar{L})= 
\bar{\phi}_{t}(\tilde{\pi}^{-1}(L))=\widetilde{\pi}^{-1}(\phi_{t}(L))
=\widetilde{\pi}^{-1}(L_{t}).$$
It follows that $$\bar{L}\cap\bar{L}_{1}\, =\, \tilde{\pi}^{-1}(L\cap L_{1}),$$ as
claimed.

\hfill $\diamond$

\subsection{The action $1$-form}

Let $L\hookrightarrow\tast$ be closed exact Lagrangian  and  let  $(\phi_{t})$ be  a
symplectic isotopy as above. Denote by $u\in H^{1}(M;\real)$ the class
$Cal(\phi_{t})$.
 We suppose that $\piu(\tast,L)=0$, using \ref{epi}.

 For $L_{t}=\phi_{t}(L)$, denote $$\Omega(L_{0},L_{1}) \,
=\, \left\{z\in\calc^{\infty}([0,1],\tast)\, |\, z(i)\in L_{i}, \, i=0,1\right\}.$$
  Define a $1$-form on $\omi$ by :
$$\nu_{z}(V)\, =\, \int_{0}^{1}\omega_{M}(z'(t),V(t))dt.$$
The zeroes of $\nu$  are the constant paths corresponding to the intersection points
$L_{0}\cap L_{1}$.

Let $\gamma:{\bf S}^{1}\ri\omi$ be a loop. We see
this loop as a map $\gamma:{\bf S}^{1}\times [0,1]\ri \tast$. We prove 
\begin{lem}\label{loop}
 We have : $$ \int_{\gamma}\nu\, =\, -u(\gamma({\bf
S}^{1}\times\{0\})).$$
In particular $\nu$ is closed. 
\end{lem}

\noindent\underline{Proof}

Let us evaluate $\int_{\gamma}\nu$.  We denote by $(s,t)$ the
coordinates on ${\bf S}^{1}\times [0,1]\ri \tast$. 
$$ \int_{\gamma}\nu\, =\, \int_{0}^{1}\nu(\partial\gamma/\partial s)ds\, =\,
\int_{0}^{1}\int_{0}^{1}\omega(\partial \gamma/\partial t, \partial\gamma/\partial s)
dt\,  ds \, =\, -\int_{\gamma({\bf S}^{1}\times[0,1])}\omega_{M}.$$
Using Stokes we find that $$ \int_{\gamma}\nu\, =\,
-\int_{\gamma(\bfs^{1}\times\{1\})}\lambda_{M}+
\int_{\gamma(\bfs^{1}\times\{0\})}\lambda_{M}.$$
Since $L=L_{0}$ is exact, the second integral in the right term above vanishes. The
first one writes :
$$\int_{\gamma(\bfs^{1}\times\{1\})}\lambda_{M}= 
\int_{\phi_{1}^{-1}(\gamma(\bfs^{1}\times\{1\}))}\phi_{1}^{\ast}\lambda_{M}=$$
$$=
\int_{\phi_{1}^{-1}(\gamma(\bfs^{1}\times\{1\}))}\lambda_{M}+
\int_{\phi_{1}^{-1}(\gamma(\bfs^{1}\times\{1\}))}(\phi_{1}^{\ast}\lambda_{M}-\lambda_{M}
).$$
  As above, the first integral in the right term is zero. Since 
$\phi_{1}^{\ast}\lambda_{M}-\lambda_{M}$ is a closed $1$-form belonging to the
cohomology class $u=Cal(\phi_{t})$ the second integral equals $u
(\phi_{1}^{-1}(\gamma(\bfs^{1}\times\{1\}))=
u(\gamma(\bfs^{1}\times\{0\}))$. Finally  
$$ \int_{\gamma}\nu\, =\, -u(\gamma({\bf
S}^{1}\times\{0\}))$$
and the proof of \ref{loop} is finished. 

\hfill $\diamond$

  \vspace{.2in}  

   Now let $(\bar{L}_{t})$ be the lifting of $(L_{t})$ to $\tastilde$, 
like in the
preceeding section and define $\Omega(\bar{L}_{0},\bar{L}_{1})$ 
as above. Also define
the $1$-form $\bar{\nu}$ on $\omibar$ in a similar way. 
The zeroes of $\bar{\nu}$
 are
therefore in bijection with the intersection points 
$\bar{L}_{0}\cap\bar{L}_{1}$. 

 For the canonical projection
$\pi^{\Omega}:\omibar\ri\omi$ we obviously have :
$$(\pi^{\Omega})^{\ast}\nu\, =\, \bar{\nu}.$$
Also remark that $\bar{\nu}$ is exact by \ref{loop}. Denote by $\cala$ a primitive of
$\bar{\nu}$. 

There is an action of $\piu(M)$ on $\omibar$, coming from the action of $\piu(M)$ on
$\tastilde$. We show :

\begin{lem}\label{equivariant} Let $z\in \omibar$ and  $g\in\piu(M)$. Then
$$\cala(g\cdot z)\, =\, \cala(z)-u(g).$$
\end{lem}

\noindent\underline{Proof}

Let $\bar{\gamma}$ be a path between $z$ and $g\cdot z$ in $\omibar$ (which does exist
since $\tastilde$ is simply connected). Denote  $\gamma=\pi^{\Omega}(\bar{\gamma})$. 
 We have :
$$\cala(g\cdot z)-\cala(z)=\int_{\bar{\gamma}}\bar{\nu}=\int_{\gamma}\nu.$$
Now $\gamma$ is a loop in $\omi$ which has the property that
$\gamma(\bfs^{1}\times\{t\})$ represents $g\in\piu(M)$ for every $t\in[0,1]$. 
By applying \ref{loop}, we get the desired relation.

\hfill $\diamond$

\vspace{.2in}

Alternatively, one may consider another action $1$-form. If the isotopy $(\phi_{t})$
is defined by the symplectic dual $X_{t}$ of  $\alpha+dH_{t}$ (defined by
$\omega_{M}(\cdot,X_{t})=(\alpha+dH_{t})(\cdot)$), we define a $1$-form on $\omil$
by :
$$\hat{\nu}_{z}(V)=\int_{0}^{1}\omega_{M}(z'(t),V(t))+(\alpha+dH_{t})(V(t))dt.$$
The zeroes of $\hat{\nu}$ are the flow trajectories starting from $L$ and ending in
$L$ at time $t=1$, which means that there are in bijection with $L_{0}\cap L_{1}$. 

As in \ref{loop}, for  a loop $\gamma$ in $\omil$ we have 
 $$ \int_{\gamma}\hat{\nu}\, =\, -
\int_{\bfs^{1}\times[0,1]}\gamma^{\ast}\omega_{M}+\int_{\bfs^{1}}\int_{0}^{1}(\alpha+d
H_{t})(\partial\gamma/\partial s)dtds$$

The first integral in the right term is zero, as in the proof of \ref{loop}. The
second one equals 
$$\int_{0}^{1}\int_{\bfs^{1}}\alpha(\partial\gamma/\partial s) ds
dt=\int_{0}^{1}\int_{\gamma(\cdot,t)}\alpha dt = \int_{\gamma(\cdot,0)}\alpha
 = u(\gamma(\cdot,0)),$$ since $\int_{\gamma(\cdot,t)}\alpha$ does not depend on $t$. 

This relation implies that $\hat{\nu}$ is closed. Then, as above, one defines a 
$1$-form on $\omibarl$ in the similar way and obtains that $\pi_{\Omega}^{\ast}\hat{\nu}\,
=\, d\hat{\cala}$. As in \ref{equivariant} we infer that 
any primitive $\hat{\cala}$ satisfies the relation
$$\hat{\cala}(g\cdot z)\, =\, \hat{\cala}(z)+u(g).$$  

These two approaches are strongly related, as it can be seen from the following
remark.  
\begin{rema}\label{correspondance} Let $\nu_{-}$ be the $1$-form defined in the same way
as $\nu$ on $\Omega(L_{0},\phi^{-1}_{1}(L_{0}))$. The map
$\Gamma(z)=\phi_{t}^{-1}z$ is obviously  a bijection between $\omil$ and 
$\Omega(L_{0},\phi^{-1}_{1}(L_{0}))$. We have the relation
$$\Gamma^{\ast}\nu_{-}\, =\, \hat{\nu}.$$
  The fact that this relation is given by a bijection which is defined using the form
$\nu_{-}$ and the 
isotopy $(\phi_{t}^{-1})$ (whose Calabi invariant is $-u$) explains the change of
sign in the analogue of \ref{equivariant} above. 
The proof is straightforward.\end{rema}

\subsection{The gradient}

 Let $(J_{t})_{t\in[0,1]}$ be a family of almost complex structures on $\tast$ which
are compatible with $\omega_{M}$. This means that $g_{t}(X,Y)=\omega_{M}(X, J_{t}Y)$
are Riemanian metrics on $\tast$. Define then a  metric $g^{\Omega}$ on
$\omi$ by :
$$g^{\Omega}(V,W)\, =\, \int_{0}^{1}g_{t}(V(t),W(t))dt.$$
The gradient of $\nu$ with respect to this metric is given by
$$grad^{g^{\Omega}}_{z}\nu\, =\,
J_{t}(z)\, z'.$$
The trajectories of the time dependent vector field ${\cal
X}_{t}=-grad^{g^{\Omega}}_{z}\nu$
 on $\omi$ can be seen as
maps $v$ of two variables $(s,t)$ satisfying the Cauchy-Riemann equation. More
precisely $v$ is a solution of $$(\ast)\gol\left\{\begin{array}{l}
\frac{\partial v}{\partial s}+J_{t}(v)\frac{\partial v}{\partial t}\, =\, 0\\
\\
v(s,0)\in L_{0}\\
\\
v(s,1)\in L_{1}\end{array}\right.$$

  For $v:\real\times[0,1]\ri\tast$,
 solution of $(\ast)$ one defines the energy $E(v)$ by the
formula $$E(v)\, =\, \int_{\real\times[0,1]}\left|\left|\frac{\partial v}{\partial
s}\right|\right|^{2}_{g_{t}}\, dsdt.$$
One can easily see that 
$$E(v) \, =\, \int_{\real\times[0,1]}\left|\left|\frac{\partial v}{\partial
t}\right|\right|^{2}_{g_{t}}\, dsdt\, =
\, \int_{\real\times[0,1]}v^{\ast}\omega_{M}.$$

Denote by $\calm$ the space of solutions of finite energy :
$$\calm(L_{0},L_{1})\, =\, \{v\in\calc^{\infty}(\real\times [0,1],\tast)\, |\, v\, \, 
\mbox{satisfies}\, \, (\ast)\, ;\, E(v)<+ \infty\}$$
For $x,y\in L_{0}\cap L_{1}$ define 
$$\calm(x,y) = \left\{v\in\calc^{\infty}(\real\times [0,1],\tast)\, \left|\, v\,  
\mbox{satisfies} \, (\ast)\, ;\, \begin{array}{c}\lim_{s\ri-\infty}v(s,\cdot)=x\\ 
\lim_{s\ri+\infty}v(s,\cdot)=y\end{array}\right. \right\}.$$

As in \cite{Flo} (Prop. 1.b), \cite{Oh} (Prop. 3.2), we have 
\begin{theo}\label{limit} $$\calm(L_{0}, L_{1})\, =\, \bigcup_{x,y\in L_{0}\cap 
L_{1}}\calm(x,y).$$
\end{theo}

   Note that the solutions of $\calm(L_{0},L_{1})$ with vanishing energy are
exactly the critical points of $\cala_{H}$ i.e. the constant paths given by the
intersection points $L_{0}\cap L_{1}$.
 However, for $x\in Crit(\cala_{H})$ the space $\calm(x,x)$ may
also contain solutions of non zero energy. It is useful to consider the spaces of
solutions with non vanishing energy : 

$$\calm^{\ast}(x,y)\, =\, \left\{\begin{array}{ccc}\calm(x,y)&\mbox{for}&x\neq y\\
\\
\calm(x,x)\setminus\{x\}&\mbox{for}&x=y\end{array}\right.$$

The family $(J_{t})_{t\in[0,1]}$ can be lifted to a family of compatible almost complex
 structures $(\bar{J}_{t})_{t\in[0,1]}$ on $\tastilde$. We define
then the spaces of solutions $\calm(\bar{L}_{0},\bar{L}_{1})$
and  
$\calm^{\ast}(\bar{x},\bar{y})$ 
 as above. Note that $\calm(\bar{x},\bar{x})=\{x\}$ (see the formula for the
energy below). 
The analogue of \ref{limit} remains
valid. 

Obviously, the projection  $\tilde{\pi}$ maps 
${\calm}(\bar{L}_{0},\bar{L}_{1})$ to $\calm(L_{0},L_{1})$ and
if $\tilde{\pi}(\bar{v})=v$ then $E(\bar{v})=E(v)$. Moreover, one
can easily see that  any
solution $v\in  \calm(L_{0},L_{1})$ lifts to a solution
 $\bar{v}\in {\calm}(\bar{L}_{0},\bar{L}_{1})$. 
Also remark that, if 
$\bar{v}\in {\calm}(\bar{x},\bar{y})$, then we have :
$$E(v)\, =\, -\int_{-\infty}^{+\infty}\frac{\partial}{\partial
s}\cala(\bar{v}(s,\cdot))ds\, =\, \cala(\bar{x})-\cala(\bar{y}).$$

Alternatively, we can consider the gradient of the $1$-form $\hat{\nu}$ with respect to
a family of metrics $\hat{g}^{\Omega}$, defined by a family of compatible almost
complex structures $(\hat{J}_{t})$,  as above. This approach leads to the perturbed
Cauchy-Riemann equation : 
 
$$(\ast\ast)\gol\left\{\begin{array}{l}
\frac{\partial \hat{v}}{\partial s}+
J_{t}(\hat{v})(\frac{\partial \hat{v}}{\partial t}-X_{t}(\hat{v}))\, =\, 0\\
\\
\hat{v}(s,0)\in L\\
\\
\hat{v}(s,1)\in L\end{array}\right.$$
  We define the energy of the solution $\hat{v}$ of $(\ast\ast)$ by the same formula
and then analogously the spaces of solutions $\hat{\calm}(L,L)$ and
$\hat{\calm}(x,y)$. The relation \ref{limit} remains true in this setting.  

  Again, for a solution $\bar{v}$
 of $(\ast\ast)$ joining $\bar{x}, \, \bar{y}\, \in\, \omibarl$
upstairs, we have $$E(\bar{v})\, =\, \hat{\cala}(\bar{x})-\hat{\cala}(\bar{y})$$
and these solutions are liftings of the solutions of $(\ast\ast)$ on $\tast$. 

\begin{rema}\label{newcorrespondence}  Let $\hat{v}\in \hat{\calm}(L,L)$. Define
$v(s,t)=\phi_{t}^{-1}(\hat{v}(s,t))$. Consider two families of compatible almost
complex structures $(J_{t})$ and $(\hat{J}_{t})$ on $\tast$ which are related by 
$$ \hat{J}_{t}(\cdot)\, =\, (\phi_{t})_{\ast}J_{t}(\phi_{t}^{-1})_{\ast}(\cdot).$$
Then one easily checks the relation :
$$\frac{\partial v}{\partial s} +J_{t}\frac{\partial v}{\partial t}\, =\,
(\phi_{t}^{-1})_{\ast}\left[\frac{\partial \hat{v}}{\partial s}+\hat{J}_{t}\left(
\frac{\partial \hat{v}}{\partial t} -X_{t}(\hat{v})\right)\right],$$
and, furthermore, the equality 
 $$E(\hat{v})=E(v).$$
We infer that $\hat{v}\mapsto v$ is a one-to-one correspondence between
$\hat{\calm}(L,L)$ and $\calm(L_{0},\phi_{1}^{-1}(L_{0}))$ for the given choices of the
almost complex families. 

  One can see this bijection as a direct consequence of the relation
$\Gamma^{\ast}\nu_{-}=\hat{\nu}$ from \ref{correspondance} for appropriate choices of
metrics on the paths spaces $\Omega(\cdot,\cdot)$. 
\end{rema}

\subsection{Transversality and compactness}

 Suppose that the manifolds $L_{0}$ and $L_{1}$ are transverse. Then we can prove the
following theorem in the same manner as  in \cite{Flo88} (see also 
\cite{Oh}). \

\begin{theo}\label{transversality} Under the transversality assumption above, and for a
generic choice of  $J_{t}$ the spaces 
$\calm^{\ast}(x,y)$ are finite dimensional
manifolds of local dimension $\mu(v)\, =\, $ the Maslov-Viterbo index of $v$ (see
\cite{Vit} for the definition). The same is true for 
${\calm}(\bar{x}, \bar{y})$ (note that ${\calm}
(\bar{x}, \bar{x})=\{\bar{x}\}$). The map $\tilde{\pi}$ induces an
embedding
 $$\tilde{\pi}:{\calm}(\bar{x}, \bar{y})\ri
\calm(x, y),$$
for $\tilde{\pi}(\bar{x})=x$ and $\tilde{\pi}(\bar{y})=y$. \end{theo}

\begin{rema}\label{newcorrespondence2}  Using \ref{newcorrespondence} we infer from
the preceeding theorem that the spaces 
$\hat{\calm}^{\ast}(x,y)$ are manifolds for a generic
choice of  $\hat{J}_{t}$. 
\end{rema}

The map $\sigma\mapsto v(\sigma+\cdot,\cdot)$ defines an action of $\real$ on
$\calm(x,y)$. Denote by ${\calitate}(x,y)$ the quotient 
$\calm^{\ast}(x,y)/\real$. Analogously, let
$\bar{{\calitate}}(\bar{x}, \bar{y})$ be the quotient of 
${\calm}^{\ast}(\bar{x}, \bar{y})$. As the action of $\real$ is
free on $\calm^{\ast}(x,y)$, we infer from \ref{transversality} 
 that ${\calitate}(x,y)$ and $\bar{{\calitate}}(\bar{x}, \bar{y})$ 
are finite dimensional
manifolds. 

 In order to define the differential of $C_{\bullet}(L_{0},L_{1})$, we need to study
the compactness of the trajectory spaces ${\calitate}(x,y)$. Remark that, since $\tast$ is
an exact symplectic manifold,
 there is no nonconstant holomorphic sphere $c:{\bf S}^{2}\ri\tast$.
Then, the relation $\omega_{M}|_{\pi_{2}(\tast,L_{i})} =0$ for $i=1,2$ implies that
there is no holomorphic disk $w:(D,\partial  D)\ri(\tast,L_{i})$. So bubbling does not
occure in a sequence in $\calm(L_{0},L_{1})$. In this framework, Gromov's classical
compactness result about holomorphic curves writes as follows :

\begin{theo}\label{Gromov}. For $A> 0$ denote $\calm_{A}(L_{0},L_{1})$ the space of
solutions $v\in \calm(L_{0},L_{1})$ with bounded energy $E(v)\leq A$. The space 
  $\calm_{A}(L_{0},L_{1})$ is compact in the topology $\calc^{\infty}_{loc}$. 
\end{theo}

As a consequence, we have the corresponding Floer-type compactness result which we will
use in the sequel :

 \begin{theo}\label{Floer} For $x,y\in L_{0}\cap L_{1}$ and $A> 0$ let
$$\calm_{A}^{\ast}(x,y)=\{v\in\calm^{\ast}(x,y)\, | \, E(v)\leq A\}$$ Let
$(v_{n})\subset\calm_{A}^{\ast}(x,y)$ be a sequence of solutions with constant index
$\mu(v_{n })=\mu_{0}$. Then there exist a finite collection $(z_{i})_{i=0,\ldots,k}$
of  points in $L_{0}\cap L_{1}$ with $z_{0}=x$ and $z_{k}=y$,   some solutions 
$v^{i}\in\calm_{A}^{\ast}(z_{i-1},z_{i})$ for 
$i=1,\ldots k$ and some sequences of real
numbers $(\sigma_{n}^{i})_{n}$ for $i=1,\ldots,k $ such that for all $i=1,\ldots,k$
 the sequence
$v_{n}(s+\sigma_{n}^{i},t)$ converges  towards $v^{i}(s,t)$ 
in $\calc^{\infty}_{loc}$.

  Moreover, we have the relations $$\sum_{i=1}^{k}E(v^{i})\, \leq\, A$$ and
$$\sum_{i=1}^{k}\mu(v^{i})\, =\, \mu_{0}.$$
\end{theo}

 We say in this case that (modulo a choice of a subsequence) $(v_{n})$ converges towards
the broken orbit $(v^{1},\ldots, v^{k})$.
  This theorem is proved in  \cite{Flo}, \cite{Oh} for the Hamiltonian
 case and stated in
\cite{Poz} for the non-exact one. Since there are some differences between these two
situations (for instance, the points $z_{i}$ need not to be  different here), we give a
complete proof below : \\

\noindent\underline{Proof of \ref{Floer}}

To simplify the notations, we will denote by $v_{n}(s)$ the path
$v_{n}(s,\cdot)\in\Omega(L_{0},L_{1})$.

 Denote by $d$ the distance $\calc^{0}$ on
$\Omega(L_{0},L_{1})$, i.e. $d(\alpha,\beta)=\sup_{t}\delta(\alpha(t),\beta(t))$,
where  $\delta$ is a distance associated to a fixed complete metric on $\tast$. 
Let $\epsilon>0$ be such that the balls $B(x,\epsilon)$ centered in $x\in L_{0}\cap
L_{1}$ are mutually  disjoint.  

  We may also suppose that every nonconstant 
 holomorphic strip $v\in \calm^{\ast}(x,x)$ leaves
 $B(x,\epsilon)$. This is true for $\epsilon$ small enough. Indeed, 
 the contrary would
imply that the image of $v$ belongs to a contractible neighbourhood of $x\in \tast$,
which means that it lifts to $\bar{v}\in \calm(\bar{x}, \bar{x})$ and therefore
$E(v)=0$ which is contradictory. 

Now denote by $v_{n}(\sigma_{n}^{1})$ the point 
(in $\Omega(L_{0},L_{1})$) where the orbit
$v_{n}(s)$ first leaves $B(x,\epsilon)$. It is given by :
$$\sigma_{n}^{1}\, =\, \inf \left\{s\in\real\, |\, d(x,v_{n}(s))>\epsilon\right\}.$$
According to \ref{Gromov}, there is a subsequence of $v_{n}(s+\sigma_{n}^{1})$ which
converges in $\calc^{\infty}_{loc}$ to an orbit $v^{1}\in\calm(L_{0},L_{1})$. Since
$v_{n}(\sigma_{n}^{1})\in\partial B(x,\epsilon)$ and 
for any $s<0$ we have $v_{n}(s+\sigma_{n}^{1})\in B(x,\epsilon) $, the limit satisfies
$v^{1}(s)\in \bar{B}(x,\epsilon)$ for $s<0$ and $v^{1}(0)\in\partial B(x,\epsilon)$. 
This implies that $v^{1}$ is a nonconstant solution in $\calm(x,z_{1})$ (i.e.
$v^{1}\in\calm^{\ast}(x,z_{1})$)  for some
$z_{1}\in L_{0}\cap L_{1}$. 

Let $s_{\ast}\in \real$ such that for all $s>s_{\ast}$ we have $v^{1}(s)\in
B(z_{1},\epsilon)$. It follows that for $n$ sufficiently large
$v_{n}(s_{\ast}+\sigma_{n}^{1})$ is in $B(z_{1},\epsilon)$. If the orbit $v_{n}$ does
not leave the ball $B(z_{1},\epsilon)$ for  $s>s_{\ast}+\sigma_{n}^{1}$, we infer
that $z_{1} =y$ and that $v^{1}(s)\in \bar{B}(y,\epsilon)$ for $s\geq s_{\ast}$. We
claim that the proof is finished in this case : we will establish the relations on the
energy and on the Maslov index at the end of this proof. 

  Now if $v_{n}(s+s_{\ast}+\sigma_{n}^{1})$ gets out the ball $B(z_{1},\epsilon)$ for
some $s>0$ consider the first exit point $v(\sigma_{n}^{2})$ defined by $$
\sigma_{n}^{2}\, =\,
\sup\left\{\sigma > s_{\ast}+\sigma_{n}^{1}\, |\, v_{n}(s)\in
B(z_{1},\epsilon), \, \forall s\in ]s_{\ast}+\sigma_{n}^{1}, \sigma[\, \right\}.$$
Using again \ref{Gromov} we find a convergent subsequence of
$v_{n}(s+\sigma_{n}^{2})$, whose limit is denoted by $v^{2}$. 

We want to show that the starting point of $v^{2}$ is $z_{1}$. 
Remark that $\sigma_{n}^{2}-\sigma_{n}^{1}\ri+\infty$. Indeed, if this sequence was
bounded, then on the bounded interval $[s_{\ast}, \sigma_{n}^{2}-\sigma_{n}^{1}]$ the
sequence $v_{n}(s+\sigma_{n}^{1})$ would converge uniformly towards $v^{1}$. In particular
$v_{n}(\sigma_{n}^{2})$ would be contained in the open ball $B(z_{1},\epsilon)$. On the
other hand, the definition of $\sigma_{n}^{2}$ implies that $v_{n}(\sigma_{n}^{2})\in
\partial B(z_{1},\epsilon)$, which yields a contradiction. 

Fix a number $s<0$. For $n$ sufficiently large we have 
$$\sigma_{n}^{1}+s_{\ast}\, <\, \sigma_{n}^{2}+s\, <\, \sigma_{n}^{2}.$$
According to the definition of $\sigma_{n}^{2}$, this means that
$v_{n}(s+\sigma_{n}^{2})\in B(z_{1},\epsilon)$ for all $s<0$, therefore, $v^{2}(s)\in
\bar{B}(z_{1},\epsilon)$ for all $s<0$. On the other hand, obviously $v^{2}(0)\in
\partial B(z_{1},\epsilon)$.  We infer that $v^{2}$ is a non-constant solution in
$\calm(z_{1},z_{2})$ for some $z_{2}\in L_{0}\cap L_{1}$. 

Then we argue in the same way to find limit solutions $v^{3}, v^{4}, \ldots$. Still we
have to prove that this iteration is valid only a finite number of times (which
means, as above, that $z_{k}=y$ for some $k$ and that $v_{n}(s+\sigma_{n}^{k})\in
B(y,\epsilon)$ for $s\geq s_{\ast}$). Suppose the
contrary : there exists a sequence of sequences $(\sigma^{i}_{n})$ such that
$v_{n}(s+\sigma_{n}^{i})$ converges towards $v^{i}(s)$ and moreover
$\sigma_{n}^{i+1}-\sigma_{n}^{i}$ are positive and unbounded. Fix a point
$\bar{z}_{i}\in\tastilde$ above each $z_{i}$ and lift $v^{i}$ to $\tastilde$
with $\bar{z}_{i}$ as ending point. The starting  point of the lift  
$\bar{v}^{i}$ will be
$g_{i}\bar{z}_{i-1}$ for some $g_{i}\in \piu(M)$. It follows that 
$$E(v_{i}) = E(\bar{v}_{i}) =\cala( g_{i}\bar{z}_{i-1})-\cala(\bar{z}_{i})=
\cala( \bar{z}_{i-1})-\cala(\bar{z}_{i})-u(g).$$
This energy is not zero, since $v^{i}$ is non constant, and as $u$ is rational
it takes values in a discrete set. Therefore, there is a constant $c>0$ such that
$$E(v^{i})\, >\, c\gol \forall\, i\in \bf N.$$

Let us show that $\sum_{i}E(v^{i})\leq A$ to get a contradiction. Fix an arbitrary
positive $\delta_{i}<E(v^{i})$ and real numbers $s_{i}^{-}<s_{i}^{+}$ such that
$$E(v^{i})-\delta_{i}\, =\, \int_{0}^{1}\int_{s_{i}^{-}}^{s_{i}^{+}}||\partial
v^{i}/\partial s||^{2}dsdt.$$
Due to the $\calc^{\infty}_{loc}$-convergence we have for $n$ sufficiently large
 (i.e. $n\geq
n_{0}$, where $n_{0}$ depends on $\delta_{i}$) :
$$E(v^{i})-2\delta_{i}<\int_{0}^{1}\int_{s_{i}^{-}}^{s_{i}^{+}}||\partial
v_{n}(s+\sigma_{n}^{i})/\partial
s||^{2}dsdt=\int_{0}^{1}\int_{s_{i}^{-}+\sigma_{n}^{i}}^{s_{i}^{+}+\sigma_{n}^{i}}
||\partial
v_{n}(s)/\partial s||^{2}dsdt.$$

Now, for $n$ large enough we have that
$s_{i}^{-}+\sigma_{n}^{i}>s_{i-1}^{+}+\sigma_{n}^{i-1}$ and
$s_{i}^{+}+\sigma_{n}^{i}<s_{i+1}^{-}+\sigma_{n}^{i+1}$. Summing up the preceeding
equality for all $i$ we obtain therefore $$\sum_{i}E(v^{i})-2\delta_{i}\, <\,
E(v_{n})\, \leq\, A.$$
Put $\delta_{i}=\delta/2^{i}$, to get $\sum_{i}E(v^{i})-\delta\leq A$, so, since
$\delta$ is arbitrary, $\sum_{i}E(v_{i})\leq A$, which contradicts $E(v^{i})>c$ for
all $i$. 

The required equality on the energies is obtained as above. It remains to show the
relation between the Maslov-Viterbo indices. It is a consequence of the following 

\begin{lem}\label{contractile} Let $\gamma_{n} : [-\infty,+\infty]\ri L_{0}$ the path
defined by $\gamma_{n}(s)=v_{n}(s,0)$ extended with $x$ in $s=-\infty$ and with $y$ in
$s=+\infty$. For $i=1,\ldots, k$ let $\gamma^{i}:[-\infty,+\infty]\ri L_{0}$ the
analogous paths defined by the holomorphic strips $v^{i}$. Then for $n$
large enough $\gamma_{n}$ and
$\gamma=\gamma^{1}\ast\gamma^{2}\ast\cdots\ast\gamma^{k}$ are homotopic in $L_{0}$. 

  The same is true for the paths defined on $L_{1}$ by $v_{n}$, $v^{1}, \ldots,
v^{k}$. \end{lem}

\noindent\underline{Proof} 

For $i=1,\ldots,k$ consider the sequences $(\sigma_{n}^{i})$ defined in the
 proof of \ref{Floer} above. They satisfy the following properties for some
$s_{\ast}>0$ :
\begin{itemize}
\item $v_{n}(s) \in B(x,\epsilon)$ for $s\leq \sigma_{n}^{1}$.
\item  $v_{n}(s) \in B(z_{i},\epsilon)$ for $s\in
[s_{\ast}+\sigma_{n}^{i},\sigma_{n}^{i+1}]$.
\item $v_{n}(s) \in B(y,\epsilon)$ for $s\geq s_{\ast}+\sigma_{n}^{k}$.
\item $v^{i}(s)\in B(z_{i-1},\epsilon)$ for $s\leq 0$, for all $i=1,\ldots,k$. 
\item $v^{i}(s)\in B(z_{i},\epsilon)$ for $s\geq s_{\ast}$, for all
$i=1,\ldots,k$.
\item $v_{n}(s+\sigma_{n}^{i})$ converges towards $v^{i}(s)$ uniformly on
$[0,s_{\ast}]$ (and more generally on every compact interval). 
\end{itemize}

The $\epsilon$-balls above are defined by the $\calc^{0}$-distance $d$ on
$\Omega(L_{0},L_{1})$ corresponding to a complete metric on $\tast$. 
Consider complete
metrics on $L_{0}$ and $L_{1}$, with associated distances $d_{0}$ and 
$d_{1}$. Define
the  distance $d'$ on $\Omega(L_{0},L_{1})$    by the formula $$d'(\alpha,
 \beta)\, =\,
max\{\, d(\alpha,\beta), \, d_{0}(\alpha(0),\beta(0)), \, d_{1}(\alpha(1),
\beta(1))\,
\}.$$
It is easy to see that we can write the first part of the proof of \ref{Floer} for
this new distance, so we may suppose that the properties of $v_{n}$ and $v^{i}$ above
are valid for the $\epsilon$-balls defined by $d'$. In particular, we can replace
$v_{n}$ by $\gamma_{n}$, $v^{i}$ by $\gamma^{i}$ and the distance by $d_{0}$. 

Now let $\chi:L_{0}\ri L_{0}$ be a continous map which is homotopic to the
identity and satisfies $\chi (B(z_{i},\epsilon))=z_{i}$ for all $i=1, \ldots,k$.
Obviously $\gamma'_{n}=\chi(\gamma_{n})$ is homotopic to $\gamma_{n}$ and
$\gamma'=\chi(\gamma)$ is homotopic to $\gamma$. Let us show that $\gamma'_{n}$ is
homotopic to $\gamma'$ for $n$ large enough. We know that
$\gamma_{n}(s+\sigma_{n}^{1})$ converges towards 
$\gamma^{1}(s)$ uniformly on $[0,s_{\ast}]$. Therefore, for $n$ sufficiently large,
$\gamma_{n}(\cdot+\sigma_{n}^{1})$ and $\gamma^{1}$ are homotopic in $L_{0}$. 
To see this, one has to write
$$\gamma_{n}(s+\sigma_{n}^{1})=exp_{\gamma^{1}(s)}Y_{n}(s),$$
for some $Y_{n}:[0,s_{\ast}]\ri (\gamma^{1})^{\ast}TL_{0}$ and then to consider the
homotopy $$\lambda\, \mapsto\, exp_{\gamma^{1}(s)}\lambda Y_{n}(s).$$
It follows that $\chi(\gamma_{n}(\cdot+\sigma_{n}^{1}))$ and $\chi(\gamma^{1})$ are
homotopic as paths defined on $[0,s_{\ast}]$. In the same way 
$\chi(\gamma_{n}(\cdot+\sigma_{n}^{i}))$ and $\chi(\gamma^{i})$ are homotopic for all
$i=1, \ldots, k$. Summarizing, we have for $n$ sufficiently large  :
\begin{itemize}
\item $\gamma_{n}'(s) = z_{i}$ for $s\in[s_{\ast}+\sigma_{n}^{i},\sigma_{n}^{i+1}]$
for $i=1,\ldots, k-1$ ; $\gamma'_{n}(s)=x$ for $s\leq \sigma_{n}^{1}$ ;
$\gamma'_{n}(s)=y$ for $s\geq s_{\ast}+\sigma_{n}^{k}$. 
\item
$\chi(\gamma^{i}(s))=z_{i-1}$ for $s\in [-\infty, 0]$  for $i=1,\ldots, k$ ; 
$\chi(\gamma^{i}(s))=z_{i}$ for $s\in [s_{\ast}, +\infty[$ for $i=1,\ldots, k$.
\item $\gamma'_{n}(s+\sigma_{n}^{i})$ and $\chi(\gamma^{i}(s))$ are homotopic on
$s\in[0,s_{\ast}]$ for $i=1, \ldots, k$. 
\end {itemize}

One easily infers that for an appropriate parametrisation $\gamma'_{n}$ and $\gamma'=
\chi(\gamma^{1}\ast\gamma^{2}\ast\cdots\ast\gamma^{k})$ are homotopic.
 The lemma is proved (the
argument for the paths on $L_{1}$ is completely analogous).

\hfill $\diamond$

Now we are able to finish the proof of \ref{Floer}. In \cite{Vit}, C. Viterbo
proved the following 

\begin{theo}\label{index}
Let $v, \, w\, \in \calm(x,y)$. Consider the paths $\gamma_{i}^{v}:[-\infty,
+\infty]\ri L_{i}$ and  $\gamma_{i}^{w}:[-\infty,
+\infty]\ri L_{i}$ defined by the restrictions of $v$ and $w$ to $\real\times\{i\}$,
for $i =0, 1$. Then 
$$\mu(w)-\mu(v)\, =\, \mu_{L_{1}}(\gamma^{w}_{1}\ast(\gamma_{1}^{v})^{-1})-
\mu_{L_{0}}(\gamma^{w}_{0}\ast(\gamma_{0}^{v})^{-1}),$$
where $\mu_{L_{i}}$ is the Maslov class of the Lagrangian manifold $L_{i}$ for
$i=0, 1$.\end{theo}

We apply the previous statement to the strips $v_{n}$ and $v^{1}\#\cdots \# v^{k}$. The
Maslov-Viterbo index of the latter is obviously $\sum_{i}\mu(v^{i})$. Then, by
\ref{contractile}, the loops in the right term of the relation above are null
homotopic for $n$ sufficiently large, so this term actually vanishes for large $n$.
 Therefore 
$$\mu_{0}\, =\, \sum_{i=1}^{k}\mu(v^{i}),$$ and the proof of \ref{Floer} is finished. 

\hfill $\diamond$

\subsection{The differential of the Floer complex}

Let $x, y\, \in \, L_{0}\cap L_{1}$. We define an "incidence number" $[x,y]$ like in
Novikov theory (Subsection 2.2).  We proceed as follows.
 Denote by ${\calitate}^{0}(x,y)$ the zero-dimensional component of
${\calitate}(x,y)=\calm^{\ast}(x,y)/\real$. For any $z\in L_{0}\cap L_{1}$ fix a lift
$\bar{z}\in \tastilde$. For $g\in\piu(M)$, consider the subset
${\calitate}^{0}_{g}(x,y)\subset{\calitate}^{0}(x,y)$ of trajectories which lift to
$\calitate (  $$g$$ $$\bar{x}$$ , \bar{y})$.  We need the following 

\begin{lem}\label{finite} Under the assumptions of \ref{transversality} 
for any $x,y \in L_{0}\cap L_{1}$ and $g\in\piu(M)$ the set ${\calitate}_{g}^{0}(x,y)$ is
finite. 

If $n_{g}$ denotes its cardinal modulo ${\bf Z}/2{\bf Z}$, the sum $\sum n_{g}g$
belongs to the Novikov ring $\Lambda_{-u}$. 
\end{lem}

\noindent\underline{Proof}

The elements of ${\calitate}^{0}_{g}(x,y)$ are classes of solutions $v$ which belong to the
$1$-dimensional component of $\calm^{\ast}(x,y)$. Moreover these solutions have the
same energy :
$$E(v)\, =\, \cala(\bar{x})-\cala(\bar{y})-u(g).$$
We can therefore apply \ref{Floer}. Since all the manifolds $\calm^{\ast}(x,y)$ have
dimensions greater than one (because of the free action of $\real$), a sequence of
solutions $v_{n}$ of Maslov index $\mu=1$ cannot converge towards a broken orbit.
Therefore, it admits a subsequence which converges in the sense of \ref{Floer} 
 towards a solution in
$\calm^{\ast}(x,y)$. This means that ${\calitate}^{0}_{g}(x,y)$ is compact, therefore it
has a finite number of elements. 
 
  Let us now prove that $\sum n_{g}g\, \in\, \Lambda_{-u}$. Let $C<0$. Our claim is
proved if we 
show that the set $\bigcup_{u(g)\geq C}{\calitate}^{0}_{g}(x,y)$ has finite cardinality. It
suffices to show that it is compact. A  sequence $(w_{n})$
 in this space lifts to a sequence
$(v_{n})$ in an $1$-dimensional component of
${\calm}^{\ast}(g\bar{x},\bar{y})$. 
If $u(g)\geq C$,
then the energy of the solutions $v$ above satisfy
$$E(v_{n})\, \leq\, \cala(\bar{x})-\cala(\bar{y})-C.$$
The sequence $(w_{n})$ 
 is therefore contained in $\calm^{\ast}_{A}(x,y)$ where $A>0$ is given by the right
term above. As above, none of its subsequences converges towards a splitting orbit.
Therefore, by \ref{Floer} $(v_{n})$  admits a subsequence which is convergent in
the considered $1$-dimensional component of 
$\calm^{\ast}_{A}(x,y)$ (in the sense of \ref{Floer}), which means that a 
subsequence of $(w_{n})$ converges in ${\calitate}^{0}(x,y)$. Moreover, the relation
between energies in \ref{Floer} shows that the limit lies in  
$\bigcup_{u(g)\geq C}{\calitate}^{0}_{g}(x,y)$. This space is therefore compact and zero
dimensional, i.e. of finite cardinality. 

\hfill $\diamond$

We define $$[x,y]\, =\, \sum_{g\in \piu(M)} n_{g}(x,y)g,$$
where $n_{g}(x,y)=\#{\calitate}^{0}_{g}(x,y)$ as above. Then, we consider the complex
$$C_{\bullet}(L_{0},L_{1},J_{t})\, =\, \Lambda_{-u}\, <L_{0}\cap L_{1}>,$$ with
differential $$\partial x\, =\, \sum_{y\in L_{0}\cap L_{1} }[x,y]y\, =\, 
\sum_{g\in\piu(M), \, y\in L_{0}\cap L_{1} }
n_{g}(x,y)gy.$$

To show the relation $\partial^{2}=0$ we have to prove that for each $g\in \piu(M)$
and $x,z\in L_{0} \cap L_{1}$ we have 
$$\sum_{y\in L_{0}\cap L_{1}, \, \, g',g''\in \piu(M), \,
g''g'=g}n_{g'}(x,y)n_{g''}(y,z)\, =\, 0.$$
This is a straightforward consequence of 

\begin{lem}\label{bord} Let ${\calitate}^{1}_{g}(x,z)$ be the $1$-dimensional component
of ${\calitate}_{g}(x,z)$. Denote by $\bar{\calitate}_{g}^{1}(x,z)$ the union 
$${\cal L}_{g}^{1}(x,z)\, \cup\, \bigcup_{y\in L_{0}\cap L_{1},
 \,
g''g'=g}{\calitate}^{0}_{g'}(x,y)\times{\calitate}^{0}_{g''}(y,z),$$
endowed with the topology given by the convergence towards broken orbits which was
defined in \ref{Floer}. 

Then $\calitate_{g}^{1}(x,z)$ is a compact 
$1$-dimensional manifold whose boundary is 
$\bigcup_{y\in L_{0}\cap L_{1}, \,
g''g'=g}{\calitate}^{0}_{g'}(x,y)\times{\calitate}^{0}_{g''}(y,z).$
\end{lem}

\noindent\underline{Sketch of the proof}\\

   To prove the compactness, let $(w_{n})$ be a sequence in
${\calitate}_{g}^{1}(x,z)$.
It admits a lift $(v_{n})\in \calm_{\ast}(x,z)$ such that $\mu(v_{n})=2$. 
By applying \ref{Floer} we find that modulo the choice of a subsequence 
$(v_{n})$ converges either
 towards a limit
$v^{0}\in\calm_{\ast}(x,z)$ or towards a broken orbit
$(v^{1},v^{2})\in\calm_{\ast}(x,y)\times\calm_{\ast}(y,z)$ for some $y\in L_{0}\cap
L_{1}$. For $i=0,1,2$ we denote by $w^{i}$ the projections of $v^{i}$ on the
correspondent trajectory spaces $\calitate$. 
In the first case, we infer using also \ref{contractile} that (modulo the  choice
of a subsequence) $(w^{n})$ converges towards $w^{0}\in{\calitate}_{g}^{1}(x,z)$. In the
last case we obtain using again \ref{contractile} that $(w_{n})$ converges towards
$(w^{1},w^{2})\in{\calitate}^{0}_{g'}(x,y)\times{\calitate}^{0}_{g''}(y,z)$ where $g',g''\in
\piu(M)$ satisfy $g''g'=g$. 

Conversely, starting with 
 $(w^{1},w^{2})\in{\calitate}^{0}_{g'}(x,y)\times{\calitate}^{0}_{g''}(y,z)$ one may use 
the usual gluing argument \cite{Flo2} 
to get a parametrisation $\Psi:[0,1[\ri{\calitate}^{1}(x,z)$
which satisfies $\lim_{\rho\ri 1}\Psi(\rho)=(w^{1},w^{2})$. 
Using \ref{contractile} we obtain that the image of $\Psi$ is contained in the
component ${\calitate}^{1}_{g''g'}(x,z)={\calitate}^{1}_{g}(x,z)$. 

\hfill $\diamond$ 

\begin{rema}\label{echivalent}  In a similar manner one can define a complex
$C_{\bullet}(L,L)$ which is spanned by the zeroes of the $1$-form $\hat{\nu}$, i.e.
 the trajectories of the flow of $X_{t}$  starting in $L$ for $t=0$ and ending in
$L$ for $t=1$. To define the differential one has to use the solution
spaces $\hat{\calm}(x,y)$ defined by $(\ast\ast)$ in \S3.3. Recall that in the
relation \ref{equivariant} written for $\hat{\cala}$ one has to change $-u$ en
$u$. This means that the complex $C_{\bullet}(L,L)$ is defined over the Novikov
ring $\Lambda_{u}$. We will denote it by 
$C_{\bullet}(L, \phi_{t}, J_{t})$ to emphasize its dependence on the
symplectic isotopy  and on the almost complex structure. Using the correspondence
\ref{newcorrespondence} one easily infers that the $\Lambda_{u}$-complexes 
$C_{\bullet}(L,\phi_{t},\hat{J}_{t})$ and $C_{\bullet}(L_{0},\phi_{1}^{-1}(L_{0}),
J_{t})$ are
actually isomorphic.
\end{rema}

 In the next subsection  we will only consider the
complex $C_{\bullet}(L,L)=C_{\bullet}(L, \phi_{t}, J_{t})$.

\subsection{Hamiltonian invariance}

  Denote by $H_{\ast}(L,\phi_{t},J_{t})$ the homology of 
$C_{\bullet}(L, \phi_{t}, J_{t})$. Recall that 
 the symplectic isotopy $(\phi_{t})$ is supposed to be defined by
$\alpha+dH_{t}$, with $\alpha$  closed $1$-form on $M$ and $H$  compactly
supported on $\tast\times[0,1]$. We want to show that this homology does not depend on
a generic choice of the couple $(J_{t}, H_{t})$ which means that it
 only depends on $L$ and on the
cohomology class $[\alpha]$ (see the analogous result for periodic orbits
 in  \cite{Leo}, Th. 4.3) : 

\begin{theo}\label{invariance}  For generic pairs $(H_{t},J_{t})$, $(H_{t}',J_{t}')$
there is an isomorphism $$\Psi: H_{\ast}(L, \phi_{t}^{\alpha+dH_{t}}, J_{t}) \ri
H_{\ast}(L, \phi_{t}^{\alpha+dH'_{t}}, J'_{t}).$$
\end{theo}

\noindent\underline{Proof}

The proof is similar to the one in \cite{Leo}, following the standard arguments in
\cite{Flo} and \cite{FS}. We consider a family of functions $H_{s,t}:\tast\ri\real$
and a family of compatible complex structures $J_{s,t}$
which depend smoothly on $(s,t)\in\real^{2}$ and which satisfy
$(H_{s,t},J_{s,t})=(H_{t},J_{t})$ for
$s<-R$ and $(H_{s,t},J_{s,t})=(H'_{t},J'_{t})$ for $s>R$, where $R>0$ is fixed.  In
order to define $\Psi$ we consider the space $\calm_{
H_{s,t}, J_{s,t}}(L)$ defined by 
$$\left\{v:\real\times[0,1]\ri \tast\, |\left |\begin{array}{c}\, 
\frac{\partial{v}}{\partial s}+J_{t}
(\frac{\partial{v}}{\partial t}-X^{\alpha+dH_{s,t}}_{s,t})=0\\
\\ v(s,i)\in L\, \mbox{for}\,
i=1,2, \, s\in\real\end{array}\right. \, ,\, E(v)<+\infty\right\}.$$

The restrictions of an element $v$ of $\calm_{
H_{s,t}, J_{s,t}}(L)$ to $s<-R$ resp. to $s>R$ are solutions of $(\ast\ast)$ 
corresponding to the couples $(H_{t},
J_{t})$ resp. $(H'_{t},J'_{t})$. One can then infer the analogue of \ref{limit},
namely the fact that any such $v$ converges towards a zero $x$ of the action $1$-form 
$\hat{\nu}$  when $s$ tends to $-\infty$
and towards a zero  $y$ of the action $1$-form 
$\hat{\nu}'$ (corresponding to the Hamiltonian $H'_{t}$)
 when $s$ tends to $+\infty$. Therefore $\calm_{
H_{s,t}, J_{s,t}}(L)$ is the union of the spaces
$\calm_{
H_{s,t}, J_{s,t}}(x,y)$ given by : $$\left\{v:\real\times[0,1]\ri \tast\, 
\left|\begin{array}{c}\, 
\frac{\partial{v}}{\partial s}+J_{t}
(\frac{\partial{v}}{\partial t}-X^{\alpha+dH_{s,t}})=0\\
\\ v(s,i)\in L\, \mbox{for}\,
i=1,2, \, s\in\real \\ 
\\
\lim_{s\ri-\infty}v(s,t)=x(t), \,
\lim_{s\ri+\infty}v(s,t)=y(t)\end{array} \, \right.\right\}.$$
Here $x,y$ are zeroes of the action $1$-forms $\hat{\nu}$ resp. $\hat{\nu}'$.
The analogue of \ref{transversality} is valid : For a generic choice of the couple
$(H_{s,t}, J_{s,t})$ the spaces $\calm_{
H_{s,t}, J_{s,t}}(x,y)$ are manifolds of local dimension at $v$ equal to the
Maslov-Viterbo index $\mu(v)$. We will define a morphism of $\Lambda_{u}$-
complexes $$\Psi_{\bullet}:
 C_{\bullet}(L,\phi_{t}^{\alpha+dH_{t}}, J_{t})\ri 
C_{\bullet}(L,\phi_{t}^{\alpha+dH'_{t}}, J'_{t}).$$ 
On the prescribed generators of $C_{\bullet}(L,\phi_{t}^{\alpha+dH_{t}}, J_{t})$ it 
is given by the formula : 
$$\Psi_{\bullet}(x)\, =\, \sum_{g\in\piu(M),\, y}m_{g}(x,y)gy, $$ where
$m_{g}(x,y)\in {\bf Z}/2{\bf Z}$ will be defined below. For this purpose we
have to consider the zero-dimensional components $\calm^{0}_{
H_{s,t}, J_{s,t}}(x,y)$. To count the elements of $\calm^{0}$ we need a compactness
result analogous to \ref{Floer}.  We obtain  indeed as in \ref{Floer} 
that any sequence $(v_{n})$ in
$$\calm_{H_{s,t}, J_{s,t}}(x,y;A)\, =\, \{v\in \calm_{
H_{s,t}, J_{s,t}}(x,y)\, |\, E(v)\leq A\}$$
has a subsequence which converges towards a broken orbit $(v^{1}, v^{2},
\ldots,v^{k})$. Actually, if the homotopy $(H_{s,t},J_{s,t})_{s}$ is not trivial, 
only   precisely one $v^{i}$ in the limit belongs to $\calm_{H_{s,t},
J_{s,t}}(L)$ ; the preceeding orbits $v^{1},\ldots, v^{i-1}$ are in $\calm_{
H_{t}, J_{t}}(L)$ and the last ones $v^{i+1}, \ldots, v^{k}$ are in
$\calm_{H'_{t},J'_{t}}(L)$. More precisely, for $j=1, \ldots, k$, $v^{j}$ is the
limit of $v_{n}(\cdot+\sigma_{n}^{j},\cdot)$ where  $\sigma_{n}^{i}=0$ 
 and the sequences $(\sigma_{n}^{j})$
tend to $-\infty$ for $j<i$, resp. to  $+\infty$ for $j>i$. 
(Too see this one has just
to pass to the limit in the Floer equation which defines 
$\calm_{H_{s,t},J_{s,t}}(L)$.)
The  energy and  the Maslov index of the limit satisfy the same
relations as  in \ref{Floer}. 

 As a consequence, we
have \begin{lem}\label{finite2}  
For any $A>0$ the set $\calm^{0}_{
H_{s,t}, J_{s,t}}(x,y;A)$ is finite. 
\end{lem}

\noindent\underline{Proof} 

   Any sequence $(v_{n})\in\calm^{0}_{
H_{s,t}, J_{s,t}}(x,y;A)$ has a subsecquence
which converges towards a broken orbit $(v^{1},\ldots, v^{k})$ as above.  It follows
 that $k =0$, since  non-constant orbits in $\calm_{H_{t},J_{t}}(L)$ and in 
$\calm_{H'_{t},J'_{t}}(L)$ have non zero Maslov-Viterbo indices.
Therefore the space $\calm^{0}_{
H_{s,t}, J_{s,t}}(x,y;A)$ is compact and zero-dimensional, so it is finite.

\hfill $\diamond$

\vspace{.2in}

   As in the preceeding subsection, fix a lift $\bar{x}$ in $\tastilde$ of every zero $x$
of
the action $1$-form $\hat{\nu}$ and a lift $\bar{y}$ in $\tastilde$ for every zero $y$ of
$\hat{\nu}'$. Consider for $g\in\piu(M)$ and any two zeroes $x,y$
of $\hat{\nu}$ resp. $\hat{\nu}'$ the space
$$\calm_{g,s}(x,y)\subset\calm_{H_{s,t},J_{s,t}}(x,y),$$ consisting of the orbits which
lift to $\tastilde$ starting from $g\bar x$ and ending at $\bar{y}$. The following
proposition is crucial for the proof of \ref{invariance}

\begin{propo}\label{energy} The space $\calm_{g,s}(x,y)$ is contained in $\calm_{
H_{s,t}, J_{s,t}}(x,y;A)$ for some $A>0$.
\end{propo}

\noindent\underline{Proof}

 We adapt the standard argument of \cite{Flo2} as in \cite{Leo} (see also
\cite{Flo}, \cite{FS}). Let $v\in \calm_{g,s}(x,y)$. We find an upper bound for $E(v)$. 
Note that in the
inequalities below (and actually in the definition of the energy) 
the norm is defined by the compatible metric
$\omega_{M}(\cdot,J_{s,t}(\cdot))$ (it therefore depends on $(s,t)$) :
$$E(v)=\int_{\real\times[0,1]}\left|\left|\frac{\partial v}{\partial
s}\right |\right|dsdt=\int_{\real\times[0,1]}\omega_{M}(\frac{\partial v}{\partial
s},J_{s,t}\frac{\partial v}{\partial
s})dsdt=$$
$$=\int_{\real\times[0,1]}\omega_{M}(\frac{\partial v}{\partial
s},\frac{\partial v}{\partial
t}-X^{\alpha}-X^{dH_{s,t}})dsdt=$$
$$=\int_{v}\omega_{M}
-\int_{\real\times[0,1]}(\alpha+dH_{s,t})(\frac{\partial v}{\partial s})dsdt \gol(1)$$

We see $v:[-\infty, +\infty]\times [0,1]$ as a path in $\omil$ between $x$ and $y$.
Fix $z_{0}\in\omil$ and let $w$ be a fixed path in $\omil$ which joins $y$ and
$z_{0}$. Denote by $v\#w$ the concatenation of $v$ and $w$. The path $w$ lifts  to a
path in $\tastilde$ joining $\bar{y}$ and  $\bar{z}_{0}$. 
Denote by $\hat{\cala}$ and
$\hat{\cala}'$ the primitives of $(\pi^{\Omega})^{\ast}\nu$ resp. of 
$(\pi^{\Omega})^{\ast}\nu'$ which vanish in  $\bar{z}_{0}$.
We have 
$$\int_{v\#w}\hat{\nu}-\int_{w}\hat{\nu}'=(\hat{\cala}
(\bar{z}_{0})-\hat{\cala}(g\bar{x}))-(\hat{\cala}'(\bar{z}_{0})-
\hat{\cala}'(\bar{y}))=-\hat{\cala}(g\bar{x})+
\hat{\cala}'(\bar{y}).\gol(2)$$
  
On the other hand, as in the computation at the end of \S3.2, we have  
$$\int_{v\#w}\hat{\nu}-\int_{w}\hat{\nu}'=-\int_{v\#w}\omega_{M}+
\int_{\real\times[0,1]}
(\alpha+dH_{t})(\frac{\partial w}{\partial s}) dsdt+$$
$$+\int_{\real\times[0,1]}
(\alpha+dH_{t})(\frac{\partial v}{\partial s})
dsdt+\int_{w}\omega_{M}- \int_{\real\times[0,1]}
(\alpha+dH'_{t})(\frac{\partial w}{\partial s}) dsdt=$$
$$=-\int_{v}\omega_{M}+\int_{\real\times[0,1]}\frac{\partial}{\partial
s}(H_{t}(w)+H_{t}(v)-H'_{t}(w))dsdt +\int_{\real\times[0,1]}\alpha(\frac{\partial
v}{\partial
s})dsdt=$$
$$=-\int_{v}\omega_{M}+\int_{[0,1]}H_{t}(z_{0})-H_{t}(x)-
H'_{t}(z_{0}) +H'_{t}(y)dt+\int_{\real\times[0,1]}\alpha(\frac{\partial
v}{\partial s})dsdt$$
 Denote by $C$ the term   $\int_{[0,1]}H_{t}(z_{0})-H'_{t}(z_{0})dt$ which does not
depend on $v$. The two
relations above imply : $$\hat{\cala}(g\bar{x})-
\hat{\cala}'(\bar{y})=-C+\int_{v}\omega_{M}
-\int_{\real\times[0,1]}\alpha(\frac{\partial
v}{\partial s})dsdt -\int_{[0,1]}-H_{t}(x)+H'_{t}(y)dt=$$
$$=-C+\int_{v}\omega_{M} -\int_{\real\times[0,1]}\alpha(\frac{\partial
v}{\partial s})dsdt-\int_{\real\times[0,1]}\frac{\partial}{\partial s}H_{s,t}(v)dsdt=$$
$$
=-C+\int_{v}\omega_{M} -\int_{\real\times[0,1]}(\alpha+dH_{s,t})(\frac{\partial
v}{\partial s})dsdt-\int_{\real\times[0,1]}\frac{\partial H}{\partial
s}(s,t,v)dsdt=$$
$$= -C+E(v)-\int_{\real\times[0,1]}\frac{\partial H}{\partial
s}(s,t,v)dsdt,\gol
(3)$$
using the relation $(1)$.

  The relation $(3)$ implies :
$$E(v)=C+\hat{\cala}(g\bar{x})-\hat{\cala}'(\bar{y})
+\int_{\real\times[0,1]}\frac{\partial H}{\partial s}(s,t,v)dsdt.$$
Since $\partial H/\partial s :\real\times[0,1]\times\tast\ri \real$ has compact
support, we infer that 
$$E(v)\leq \hat{\cala}(g\bar{x})-\hat{\cala}'(\bar{y})+K=
\hat{\cala}(\bar{x})-\hat{\cala}'(\bar{y})+u(g)+K \gol(4),$$ for some  
$K$ which does not
depend on $v$. It follows that  $\calm_{g,s}(x,y)$ is contained in $\calm_{
H_{s,t}, J_{s,t}}(x,y;A)$ for some $A>0$, as required. 

\hfill $\diamond$

\vspace{.2in}

A straightforward consequence of \ref{energy} is that the set
$$\calm_{g,s}^{0}(x,y)=\calm_{g,s}(x,y)\cap\calm^{0}_{H_{s,t},J_{s,t}}(x,y)$$
has finite cardinality. This enables us to define the morphism $$\Psi_{\bullet}:
 C_{\bullet}(L,\phi_{t}^{\alpha+dH_{t}}, J_{t})\ri 
C_{\bullet}(L,\phi_{t}^{\alpha+dH'_{t}}, J'_{t}).$$ by the formula
$$\Psi_{\bullet}(x)\, =\, \sum_{g\in\piu(M),\, y}m_{g}(x,y)gy, $$
where $m_{g}(x,y)$ is the parity of $\calm_{g,s}^{0}(x,y)$. Note that by
\ref{energy} the coefficients $\sum_{g\in\piu(M)} m_{g}(x,y)g$ belong to
$\Lambda_{u}$. Indeed, for any $B\in\real$, the relation $(4)$ above shows that
$\bigcup_{u(g)<B}\calm_{g,s}^{0}(x,y)$ is contained in
$\calm^{0}_{H_{s,t},J_{s,t}}(x,y;A)$ for some positive constant $A$,
 so it is finite, according to \ref{finite}.

  The fact  that $\Psi_{\bullet}$ commutes with the differentials can be 
proved in the usual way, by studying the compactness  of the $1$-dimensional
components of  $\calm_{H_{s,t},J_{s,t}}(x,y)$ like in the proof of \ref{Floer}. 
(a sequence in this space either admits a convergent subsequence, 
or converges towards a broken orbit
$(v^{1},v^{2})$). The proof is similar to \ref{bord}. 

  Finally, to show that $\Psi_{\bullet}$ induces an isomorphism in homology, one again
uses the standard method of Floer theory \ref{Floer} (construct a morphism 
$\Gamma_{\bullet}:
 C_{\bullet}(L,\phi_{t}^{\alpha+dH'_{t}}, J'_{t})\ri 
C_{\bullet}(L,\phi_{t}^{\alpha+dH_{t}}, J_{t})$ and than prove that
$\Psi_{\bullet}\Gamma_{\bullet}$ and $\Gamma_{\bullet}\Psi_{\bullet}$ are homotopic to
the identity, using a two-parameter homotopy $H_{r,s,t}$). 

The proof of \ref{invariance} is now finished. 

\hfill $\diamond$

\begin{rema}\label{notyet} The complex $C_{\bullet}(L, \phi_{t}, J_{t})$ defined in the previous
subsection is free over the Novikov ring $\Lambda_{-u}$. As we remarked in \ref{echivalent} one can infer
the existence of a similar complex over $\Lambda_{u}$. We showed above that the homology of these
complexes only depends on $L$ and on $u$. But the goal of our theorem \ref{conclusion} was
more general, namely the existence of a complex which is free over $\Lambda_{p^{\ast}u}$ (spanned by the
intersection points $L\cap\phi_{t}(L)$). In the next subsection we show how to adapt the previous
construction in order to get this conclusion.
\end{rema}

\subsection{The Floer-Novikov complex over $\Lambda_{p^{\ast}u}$.}

The idea is the following. Consider the intersection points $L\cap \phi_{1}(L)$, 
viewed as points in $L$.
For two such points $x, y$, any holomorphic strip $v\in \calm(x,y)$ defines a path 
$\gamma:]-\infty,
+\infty[ \ri L$ which joins $x$ and $y$ : 
$$\gamma(s) \, =\, v(s,0).$$
 Look at the collection of intersection points and take the paths $\gamma$ as 
above, defined by the
strips $v$ which belong to the one-dimensional component of $\calm(x,y)$ 
(which correspons to the
zero dimensional component of $\calitate(x,y)$). This collection of points
 and paths joining them 
is sufficient to re-construct the complex 
$C_{\bullet}(L, \phi_{t}, J_{t})$ : one just has to fix lifts $\bar{x}\in\bar{L}$ for any point and then
lift the lines $\gamma$ from $g\bar{x}$ to $\bar{y}$. We get thus the same "incidence number"
$[x,y]\in\Lambda_{-u}$ as
above. Therefore we obtain the same complex. 

  Now, instead of lifting these lines to the covering space $\bar{L}$, 
we lift them to the universal cover
$\widetilde{L}$. If we start with fixed lifts $\tilde{x}\in \widetilde{L}$ 
of the intersection points, we
get thus a new incidence number $[x,y]^{\sim}$ which belongs to the Novikov ring 
$\Lambda_{-p^{\ast}u}$. 
This enables us to define the desired complex. \\

This idea can be formalized as follows :    Let $L$ be a closed manifold and let 
$\calc$ be a finite set of points on $L$.
Consider a (possibly infinite) collection $\cal{G}$ of paths 
$\gamma: [-\infty, +\infty]\ri L$ such that
$\gamma(\pm \infty)\in\calc$. Let $p:\piu(L)\ri G$ be an epimorphism 
onto a group $G$ and let $u:G\ri{\bf
Z}$ be a group morphism. Consider the covering space $\bar{L}\ri L$ 
associated to $Ker(p)\subset \piu(L)$.
For any $x\in \calc$ fix a lift $\bar{x}\in\bar{L}$. For $x,y\in\calc$, and $g\in G$,
denote $\calitate^{0}_{g}(x,y)$ the set of the paths
 $\gamma\in{\cal G}$ which have a lift in $\bar{L}$ which
joins $g\bar{x}$ and $\bar{y}$. Denote by $\Lambda_{u}$ the Novikov ring 
${\bf Z}/2\, [G]_{u}$. Now we prove :

\begin{propo}\label{lift}
Suppose that :\\
a) For any $x,y\in \calc$ and $g\in G$ the space $\calitate^{0}_{g}(x,y)$ is finite and  
$$[x,y]=\sum_{g\in G}\#_{2}\calitate^{0}_{g}(x,y)g\, \in\, \Lambda_{-u}.$$
b) The formula $$\partial x =\sum_{y\in \calc}[x,y]y$$
defines the differential of a $\Lambda_{-u}$-free complex $C_{\bullet}$ 
spanned by $\calc$. Note that this is
equivalent to the fact that for any $x,z\in\calc$ and $g\in G$ the space 
$$\bigcup_{y\in \calc,
 \,
g''g'=g}{\calitate}^{0}_{g'}(x,y)\times{\calitate}^{0}_{g''}(y,z)$$ has 
an even number of elements.
\\
c)  The space above is a disjoint union of sets with two elements  
$\{(\gamma_{1},\gamma_{2}),(\gamma'_{1},\gamma'_{2})\}$ with the property 
that the paths 
 $\gamma_{1}\ast\gamma_{2}$
and $\gamma_{1}'\ast\gamma_{2}'$ are homotopic in $L$. 

Then there exists a free $\Lambda_{-p^{\ast}u}$-complex $\widetilde{C}_{\bullet}$,
 spanned by $\calc$ such that $p:\piu(L)\ri G$ induces a morphism from 
$\widetilde{C}_{\bullet}$ to
$C_{\bullet}$ via the natural ring morphism $p:\Lambda_{-p^{\ast}u} \ri \Lambda_{-u}$.
\end{propo}

\noindent\underline{Proof}

Fix  lifts $\tilde{x}\in \widetilde{L}$ of the points $x\in\calc$. For any 
$h \in \piu(L)$, denote by $\widetilde{\calitate}^{0}_{h}(x,y)$ the set of paths 
$\gamma\in{\cal G}$ which lift in $\widetilde{L}$ to  paths joining  
$h\tilde{x}$ and $\tilde{y}$. Then define 
$$[x,y]^{\sim}=\sum_{h\in \piu(L)}\#_{2}\widetilde{\calitate}^{0}_{h}(x,y)h.$$ 
Remark that $[x,y]^{\sim}$ belongs to $\Lambda_{-p^{\ast}u}$.
 Indeed it is obvious that for any $g\in G$ we have 
$$\calitate^{0}_{g}(x,y)\, =\, 
\bigcup_{h\in \piu(L), \, p(h)=g}\widetilde{\calitate}^{0}_{h}(x,y).$$
 Now define 
$$\partial x =\sum_{y\in \calc}[x,y]^{\sim}y.$$
Like in the case of  the statement b in \ref{lift} above, 
the relation $\partial\circ\partial=0$ 
is equivalent to the fact that  for any $x,z\in \calc$ and $h\in \piu(L)$ 
the number of elements of the set 
$$\bigcup_{y\in \calc,
 \,
h''h'=h}{\widetilde{\calitate}}^{0}_{h'}(x,y)
\times{\widetilde{\calitate}}^{0}_{h''}(y,z)$$ 
is even. Let $(\gamma_{1},\gamma_{2})$ be an element of this set. 
Then the set equality above implies that 
$$(\gamma_{1},\gamma_{2})\, \in\, \bigcup_{y\in \calc,
 \,
g''g'=g}{\calitate}^{0}_{g'}(x,y)\times{\calitate}^{0}_{g''}(y,z),$$ 
where $g=p(h)$, $g'=p(h')$ and $g''=p(h'')$. Let $(\gamma_{1}',\gamma_{2}')$ 
as in the hypothesis c of \ref{lift}. 
Since $\gamma_{1}\ast\gamma_{2}$ and $\gamma_{1}'\ast\gamma_{2}'$ are homotopic, 
it follows that $$\gamma_{1}'\ast\gamma_{2}'\, \in\, \bigcup_{y\in \calc,
 \,
h''h'=h}{\calitate}^{0}_{h'}(x,y)\times{\calitate}^{0}_{h''}(y,z).$$
Therefore, using c, the latter has an even number of elements. 
This implies that $\partial\circ\partial=0$ and proves the proposition. 

\hfill $\diamond$

We apply \ref{lift} to $\calc=L\cap\phi_{1}(L)$ and 
$\cal G$ defined by the paths $\gamma(s)=v(s,0)$, where $v$ 
are the holomorphic strips belonging to the zero-dimensional components 
of the trajectory spaces $\calitate(x,y)$. In order to check the hypothesis c 
of \ref{lift} recall that the space 
$$\bigcup_{y\in \calc,
 \,
g''g'=g}{\calitate}^{0}_{g'}(x,y)\times{\calitate}^{0}_{g''}(y,z),$$ is the 
boundary of a one-dimensional closed manifold (as it was shown in \ref{bord}),
 so its elements can be viewed as the disjoined union of the boundaries 
(consisting of two elements) of the connected components of this manifold. 
Such a couple $(\gamma_{1},\gamma_{2}), (\gamma_{1}',\gamma_{2}')$ has 
the property that $\gamma_{1}\ast\gamma_{2}$ and $\gamma_{1}'\ast\gamma_{2}'$ 
are homotopic in $L$. This is an immediate consequence of \ref{contractile}. 
We get thus a free $\Lambda_{-p^{\ast}u}$-complex, as in \ref{lift}. 

Using the equivalent approach \ref{echivalent}, we obtain a  
free $\Lambda_{p^{\ast}u}$-complex 
$\widetilde{C}_{\bullet}(L, \phi_{t},J_{t})$ spanned by the intersection points 
$L\cap\phi_{1}(L)$, as claimed in Theorem \ref{conclusion}. 
 In order to show that its homology only depends on $L$ and on $u$ one has to 
prove that the the morphism of $\Lambda_{u}$-complexes 
$$\Psi_{\bullet}:
 C_{\bullet}(L,\phi_{t}^{\alpha+dH_{t}}, J_{t})\ri 
C_{\bullet}(L,\phi_{t}^{\alpha+dH'_{t}}, J'_{t})$$ defined in the previous 
subsection lifts to a morphism between the corresponding 
 $\Lambda_{p^{\ast}u}$-complexes. 
This is obtained using the same argument as above. We also prove that the 
lifted morphism yields an isomorphism in homology in an analogous way.

  The goal of this section, Theorem \ref{conclusion}, is now achieved.

\section{Floer homology and Novikov homology}

   Denote by $FH(L, u)$ the homology of the Floer complex
$\widetilde{C}_{\bullet}(L, \phi_{t}, J_{t})$ defined in Subsection 3.7. 
Denote by $H(L,p^{\ast}u)$   the Novikov homology of $L$ and of the class $p^{\ast}u$, 
 where $p:L\ri M$ is the projection on the base space of $\tast$.
The aim of this section is to show that

\begin{theo}\label{iso} $FH(L,u)$ is isomorphic to the Novikov homology
$H(L,p^{\ast}u)$.\end{theo}

Our 
results \ref{main} and \ref{coro} will be inferred from this theorem.

\subsection{Proof of \ref{iso}}

  Again, we follow the ideas of \cite{Leo}. We prove

\begin{propo}\label{epsilon} Let $u\in H^{1}(M)$   There exists an $\epsilon >0$
(depending on $u$) such that :\\
a) For all $|\sigma|<\epsilon$, 
$$FH(L, (1+\sigma)u)\, \approx\, FH(L,u).$$
b) $FH(L,\sigma u)\, \approx\, H(L,p^{\ast}u)$
\end{propo}

Proposition \ref{epsilon} immediately implies \ref{iso} since the set
$$E=\{\, \sigma\in]0,+\infty[\, |\, FH(L, \sigma u)\approx H_{\ast}(L,p^{\ast}u)\}$$
is non empty, open and with open complementary, so it equals $]0,+\infty[$. \\

\noindent\underline{Proof of \ref{epsilon}}

Consider the $\Lambda_{u}$-complex $C_{\bullet}(L,\phi_{t}, J_{t})$ 
defined in Subsection 3.5 (using
\ref{echivalent}). We showed in \S3.6 that its homology only depends on 
$L$ and on $u$. 
Denote this homology by $\bar{FH}(L,u)$. Also consider the Novikov homology
 associated to
$p^{\ast}u$ and to the covering ${\bar L}\ri L$ (this covering
 was defined in Subsection 3.1 as the
pull-back to $L$ of the covering $\tastilde\ri \tast$ ; it corresponds to
$Ker(p)\subset\piu(L)$). This homology, defined as explained in Remark
\ref{novikovgeneral},  will be denoted  by $H(\bar{L}\ri L, p^{\ast}u)$. 
In order to show \ref{epsilon} we prove first the analogous result for the Floer
homology $\bar{FH}(L,u)$, namely :  

\begin{propo}\label{epsilon2} Let $u\in H^{1}(M)$   There exists an $\epsilon >0$
(depending on $u$) such that :\\
a) For all $|\sigma|<\epsilon$, 
$$\bar{FH}(L, (1+\sigma)u)\, \approx\, \bar{FH}(L,u).$$
b) $\bar{FH}(L,\sigma u)\, \approx\, H(\bar{L}\ri L,p^{\ast}u)$
\end{propo}
 
Then, using the same method as in Subsection 3.7, we show that the isomorphisms
\ref{epsilon2}.a and \ref{epsilon2}.b can be lifted to the isomorphisms 
\ref{epsilon}.a
and respectively \ref{epsilon}.b, proving thus \ref{epsilon}. \\

\noindent\underline{Proof of \ref{epsilon2}}   

 Let $\alpha\in u$ a closed $1$-form. 
In order to compute $\bar{FH}(L,u)$ we choose a
symplectic isotopy $\phi_{t}^{\alpha+dH_{t}}$, as follows. Let
$\psi_{t}:\tastl\ri\tastl$
the symplectic isotopy defined by $\psi_{t}(x)=x+tp^{\ast}\alpha$. Now we use the
following well-known result 

\begin{lem}\label{Laudenbach}
If $L\subset\tast$ is exact Lagrangian there exists a (non-proper) 
symplectic embedding
$\Phi:\tastl\ri\tast$, extending the given embedding of $L$. 
In particular, $\Phi^{\ast}\lambda_{M}-\lambda_{L}$ is an exact one
form $dG$.
\end{lem}

\noindent\underline{Proof} 

By Weinstein's theorem there is a symplectic embedding $\Phi$ of a tubular 
neighbourhood $U$ of
$0_{L}$ whose restriction to $0_{L}$ is the given embedding of $L$. Since $L$ is exact
the difference $\lambda_{M}-(\Phi^{-1})^{\ast}\lambda_{L}$ is an exact $1$-form on 
 $\Phi(U)$. This enables one to extend $(\Phi^{-1})^{\ast}\lambda_{L}$ to a primitive
of $\omega_{M}$ on $\tast$. The symplectic dual of this primitive is a vector field
whose restriction on $\Phi(U)$ is the image of the canonical Liouville vector field on
$U\subset\tastl$. Denote by $\xi_{t}$ the flow of this vector field and by $\rho_{t}$
the flow of the canonical Liouville vector  field on $\tastl$. Then 
the embedding $\Phi$ is defined by the formula $$\Phi(x)=\xi_{t}\circ\Phi|_{U}\circ
\rho_{-t}(x),$$ where $t>0$ is sufficiently large to ensure $\rho_{-t}(x)\in U$. It is
easy to see that this definition does not depend on $t$  and that $\Phi$ is an exact
symplectic embedding as claimed.

 \hfill $\diamond$

Consider now the Lagrangian isotopy $\Phi\circ\psi_{t}|_{L} :L\ri\tast$. There is a
symplectic isotopy $(\phi_{t})$ on $\tast$ which extends $\Phi\circ\psi_{t}$. To see
this, one has to consider the isotopy $\chi_{t}:\tast\ri\tast$ defined by
$\chi_{t}(x)=x+t\alpha$. It is easy to see that
$(\chi_{-t}\psi_{t})^{\ast}\lambda_{M}$ is an exact $1$-form on $L$, so
$\chi_{-t}\psi_{t}$ is
an exact Lagrangian isotopy. Consider a Hamiltonian  extension $(\Gamma_{t})$ of 
$\chi_{-t}\psi_{t}$. Then $\chi_{t}\Gamma_{t}$ is an extension of $\psi_{t}$.
Therefore we can consider a symplectic extension of $\Phi\circ\psi_{t}$, which we
denote by $(\phi_{t})$.  The Calabi invariant
of the extension  is clearly $u=[\alpha]= Cal(\psi_{t})$ since
$(\Gamma_{t})$ is Hamiltonian and $\Phi$ is an extension of an exact Lagrangian
embedding. 

  Using \ref{Calabi}, we may suppose that $(\phi_{t})$ is defined by $\alpha+dH_{t}$,
where $H:[0,1]\times\tast\ri \real $ is compactly supported. We will use this isotopy
for the definition of the Floer complex. Note that the intersection points
$L\cap\phi_{t}(L)$ are  the zeroes of $p^{\ast}\alpha$ in $L$ and therefore
they are fixed with respect to $t$. In other words, the zeroes of the action $1$-form
$\hat{\nu}$ are constant paths in $\omil$. Note also that when $\alpha$ is Morse
(which we will always suppose) the
intersections $L\cap \phi_{t}(L)$ are transverse, so the isotopy $(\phi_{t})$ is
generic in this sense. 

We will also need the following Palais-Smale-type lemma (see \cite{Leo}, Lemma 5.1)

\begin{lem}\label{Palais}
Let $L_{t}=\phi_{t}(L)$, as above and denote by $\{x_{1}, \ldots, x_{k}\}$ the
intersection points $L\cap L_{t}$ for $t>0$. Fix a ball $B_{i}\subset\tast$ around
each $x_{i}$ and denote by $B$ the union $\bigcup_{i}B_{i}$. Then there exist $c>0$
such that for any  smooth $z\in \omil$, whose image is not contained in $B$ we have 
$$||z'(t)-X^{\alpha+dH_{t}}(z(t))||_{L^{2}}\, \geq \, c.$$
\end{lem}

\noindent\underline{Proof}

  The norm $L^{2}$ in the statement above is defined using a fixed complete metric on
$\tast$. 
Suppose the contrary of \ref{Palais} : there exist a sequence 
$(z_{n})\in\omil$ of paths whose
images are not contained in $B$, such that 
$$\lim_{n\ri+\infty}||z_{n}'(t)-X^{\alpha+dH_{t}}(z_{n}(t))||_{L^{2}}\, =\, 0.$$
Since  $H_{t}$ is compactly supported  and $\alpha$ is defined on $M$ the norm
$||X^{\alpha}(z)+X^{dH_{t}}(z)||_{L^{2}}$ is bounded uniformly with respect to $z$, so
there is a constant $K>0$ such $ ||z_{n}'||_{L^{2}}\leq K$ for all $n\in \bf N$. Let
$d$ be the distance defined on $\tast$ by the metric we considered. For arbitrary
$t_{0}<t_{1}$ in $[0,1]$ we have 
$$d(z_{n}(t_{0}), z_{n}(t_{1}))\leq
\int_{t_{0}}^{t_{1}}||z_{n}'(t)||dt=\int_{0}^{1}1_{[t_{0},t_{1}]} ||z_{n}'(t)||dt\leq$$
$$\leq ||z_{n}'||_{L^{2}}||1_{[t_{0},t_{1}]}||_{L^{2}}\leq K\sqrt{t_{1}-t_{0}}.$$
The family $(z_{n})$ is therefore equicontinous. Since $(z_{n}(0))\in L$ admits a
convergent subsequence, we may apply Arzela-Ascoli to get a subsequence of $(z_{n})$
which converges towards some $z_{\infty}\in\omil$ in the topology
 $\calc^{0}([0,1],\tast)$.  It follows that $X^{\alpha+dH_{t}}((z_{n})$ converges
towards $X^{\alpha+dH_{t}}(z_{\infty})$ in the topology $\calc^{0}$ and in particular
in the norm $L^{2}$.  

But $||z_{n}'(t)-X^{\alpha+dH_{t}}(z_{n}(t))||_{L^{2}}$ converges to zero, so we have
 $$\lim_{n\ri+\infty}||z_{n}'(t)||_{L^{2}}\, =\,
||X^{\alpha+dH_{t}}(z_{\infty}(t))||.$$
Embed $\tast$ is some Euclidean space $\real^{N}$ and see   the vectors fields 
 in the equality above as
elements of   $\calc^{0}([0,1],\real^{N})$ (depending on the variable $t$). 
Obviously the last convergence
is valid in $L^{2}([0,1],\real^{N})$. Then one can write for $t\in[0,1]$  
$$z_{n}(t)-z_{n}(0)=\int_{0}^{t}z'_{n}(\tau)d\tau=\int_{0}^{t}z_{n}'(\tau)-X^{\alpha+d
H_{t}}(z_{\infty}(\tau))d\tau +\int_{0}^{t}X^{\alpha+dH_{t}}(z_{\infty}(\tau))d\tau.$$
Using the Cauchy-Schwarz inequality, we find as above 
$$ ||\int_{0}^{t}z_{n}'(\tau)-X^{\alpha+dH
_{t}}(z_{\infty}(\tau))d\tau|| \, \leq \, ||z_{n}'-X^{\alpha+d
H_{t}}(z_{\infty})||_{L^{2}}\sqrt{t},$$
in particular this integral converges to zero.  We infer that 
when $n$ goes to $+\infty$ the
preceeding equality writes  :
$$z_{\infty}(t)-z_{\infty}(0)\,
=\, \int_{0}^{t}X^{\alpha+dH_{t}}(z_{\infty}(\tau))d\tau.$$
In particular, $z_{\infty}$ is $\calc^{1}$ (and hence $\calc^{\infty}$, by an
obvious bootstrapping argument) and satisfies
$$z_{\infty}'\, =\, X^{\alpha+dH_{t}}(z_{\infty}).$$
This means that $z_{\infty}$ is a zero of the action $1$-form 
$\hat{\nu}$, hence it is a constant path which belongs to $\{x_{1},\ldots, x_{k}\}$. 
But this is contradictory, since the image of $z_{n}$ is not contained in $B$, so the
sequence $(z_{n})$ cannot converge towards an element $z_{\infty}\in \{x_{1}, \ldots,
x_{k}\}$. 

\hfill $\diamond$

Now we are able to give the \\
\\
\noindent\underline{Proof of \ref{epsilon2}.a} 

Recall that we have a generic isotopy
$(\phi_{t})$ which is defined by $X^{\alpha+dH_{t}}$. The intersection points
$L\cap\phi_{t}(L)$ are fixed with respect to $t$ ; we denoted them $\{x_{1},\ldots,
x_{k}\}$. As in Lemma \ref{Palais} we fix a collection of balls 
 around these intersection
points and we denote its union by $B$.  We also 
consider the constant $c$ given by this lemma. Recall also that $u\in
H^{1}(M)$ is the cohomology class of $\alpha$.

 Choose $\eta\in u$ such that $\eta|_{B}=0$ 
and fix $\epsilon>0$ such that $$\epsilon||\eta||<c/3.$$ Then pick
$\sigma<\epsilon$ and consider the symplectic isotopy $(\psi_{t})$ defined by
$X^{\alpha+\sigma\eta+dH_{t}}$. The constant $\epsilon>0$ is chosen small enough to 
 ensure  that
$\psi_{1}(L)$ is still transverse to $L$ (actually we may even suppose that the
intersection points are $\{x_{1}, \ldots, x_{k}\}$ but this is not needed in the
proof).    
 Now fix a
compatible almost complex structure $J$ on $\tast$ 
which yields a complete metric.
Then choose
compatible complex structures $J_{t}$ and $J'_{t}$ such that the couples
$(\alpha +dH_{t},J_{t})$ and $(\alpha +dH_{t}+\sigma\eta ,J'_{t})$ satisfy the
transversality assumption of \ref{transversality}.
By genericity, we may suppose that $||J_{t}-J||<\delta$ and $||J'_{t}-J||<\delta$ 
where $\delta>0$ is a (small) constant which will be specified later. The norm here is
defined by the metric $g_{J}$, induced by $J$.  
Like in the previous section,
define the  $\Lambda_{u}$-complexes
 $C_{\bullet}(L, \phi_{t}^{\alpha+dH_{t}}, J_{t})$  and 
$C_{\bullet}(L, \psi_{t}^{\alpha+\sigma\eta+dH_{t}}, J'_{t})$ (we use here that
$\Lambda_{u}=\Lambda_{\tau u}$ for any $\tau>0$). 

To finish the proof of \ref{epsilon2}.a we have to prove that the homologies of these
two complexes are isomorphic. We proceed as in Section 3.6 by constructing a homotopy 
between the pairs $(\alpha+dH_{t}, J_{t})$ and $(\alpha
+\sigma\eta+dH_{t},J'_{t})$. Denote by $(\alpha+\chi(s)
\sigma\eta+dH_{t},J_{t,s})$ this
homotopy. Here $\chi$ is a monotone increasing smooth function on $\real$ which
vanishes for $s\leq -R$ and equals $1$ for $s\geq R$. We chose the homotopy of almost
complex structures such that :
 $ J_{s,t}=J_{t}$ for $s\leq -R$ and   
$ J_{s,t}=J'_{t}$ for $s\geq R$. We may  also suppose that for all $s\in\real$
 $||J_{s,t}-J||<\delta$.  In order to define
a morphism between the two complexes above we need to consider the solutions
$v:\real\times[0,1]\ri\tast$ of the
system
$$\left\{\begin{array}{c}
\frac{\partial v}{\partial s} +J_{s,t}\left(\frac{\partial v}{\partial
t}-X^{\alpha+dH_{t}+\chi(s)\sigma\eta}(v)\right)\, =\, 0\\
\\ 
v(s,i)\in L\, \, \mbox{for}\, \, i=0,1.\\
\\
E(v)<+\infty\end{array}\right., $$ where $E(v)$ is the energy of $v$ with respect to
the norm defined by $J$ (or equivalently, to the norm defined by $J_{s,t}$).
 As in the previous section, each solution $v$ of this
system satisfies $\lim_{s\ri-\infty}v(s,t)=x(t)$ and $\lim_{s\ri+\infty}v(s,t)=y(t)$,
where $x(t)\in\omil$ is an orbit of $X^{\alpha +dH_{t}}$ and $y(t)\in\omil$ is an
orbit of $X^{\alpha+\sigma\eta+dH_{t}}$ (equivalently, they are zeroes of the
corresponding action $1$-forms). The genericity assumptions ensure that the spaces
$\calm_{\chi,J_{s,t}}(x,y)$ of solutions with the indicated limit conditions are
manifolds of local dimension $\mu(v)$.   As in \ref{Floer}, the 
 zero dimensional subspaces
$\calm^{0}_{\chi,J_{s,t}}(x,y;A)$
of solutions   with energy uniformly bounded by $A$ are compact and those
of dimension $1$ are compact up to breaking into $(v^{1},v^{2})$ where only one of the
$v^{i}$'s is a solution of the equation above, the other being a solution of the
Floer equation corresponding either to $(\alpha+dH_{t},J_{t})$ or to 
$(\alpha+\sigma\eta+dH_{t},J'_{t})$. We want to define a morphism
$$\Gamma_{\bullet}:   C_{\bullet}(L, \phi_{t}^{\alpha+dH_{t}},
J_{t})\ri C_{\bullet}(L, \psi_{t}^{\alpha+\sigma\eta+dH_{t}}, J'_{t})$$
by the formula $$\Gamma (x)=\sum_{g\in\piu(M),y}m_{g}(x,y)gy,$$ 
where $m_{g}(x,y)$ is the number mod $2$ of elements of the space
$\calm^{0}_{g,s}(x,y)\subset 
\calm^{0}_{\chi,J_{s,t}}(x,y)$ of solutions which lift to $\tastilde$ starting form
$g\bar{x}$ and ending at $\bar{y}$ (as previously, we fixed lifts $\bar{x}$ and
$\bar{y}$ for all   the zeroes of the two action $1$-forms). 

  The crucial point is the following statement, analogous to \ref{energy}. It implies
 that
$\calm^{0}_{g,s}(x,y)$ is finite and that   for any $y$ the sum
$\sum_{g\in\piu(M)}m_{g}(x,y)g$ belongs to $\Lambda_{u}$ :

\begin{propo}\label{crucial} The space $\calm_{g,s}(x,y)$ is contained in 
$\calm_{\chi, J_{s,t}}(x,y;A)$ for some $A>0$. 
\end {propo} 

\noindent\underline{Proof}

Let $v\in\calm_{g,s}(x,y)$. As in the previous section $\hat{\cala}$ and
$\hat{\cala}'$ are the
primitives of the pull-backs to $\tastilde$  of the two action $1$-forms $\hat{\nu}$
and $\hat{\nu}'$. Denote by $\bar{v}$ the lift of $v$ to $\tastilde$.  The proof of
\ref{crucial} is implied by the estimate 
$$E(v)\, \leq\, 3[\hat{\cala}(\bar{x})-\hat{\cala}(\bar{y})+u(g)].\gol (2)$$
Let us prove this inequality. In the estimations below the scalar product $<,>$ is
$g_{J}(\cdot,\cdot)=\omega_{M}(\cdot, J\cdot)$ for the fixed structure $J$ and the
norm $||\cdot||$ is defined by this metric.
  We have
$$\hat{\cala}(g\bar{x}) -\hat{\cala}(\bar{y})= -\int_{-\infty
}^{+\infty}\frac{\partial}{\partial s}\hat{\cala}(\bar{v}(s,\cdot))ds=$$
$$= 
 -\int_{-\infty
}^{+\infty}\hat{\nu}(\frac{\partial v}{\partial s})ds   
=-\int_{-\infty}^{+\infty}\left<\frac{\partial v}{\partial s}, 
grad^{g_{J}}_{v}\hat{\nu}\right>ds
=$$
$$=-\int_{-\infty}^{+\infty}\int_{0}^{1}\left<\frac{\partial v}{\partial s},
J(\frac{\partial
v}{\partial t} -X^{\alpha+dH_{t}}(v))\right>dtds.\gol (3) $$

In order to prove $(2)$ we have to find a lower bound for the right term of the
previous equality. Let us fix $s\in \real$. Recall that $B$ is  a fixed union of balls
around the intersection points $L\cap\phi_{t}(L)$ such that Lemma \ref{Palais} is
valid.  We consider the following cases :\\
\\
$1^{\circ}$ $Im(v(s,\cdot))\subset B$.

 Using that $v$ is a solution of the
parametrized Floer equation, we get :
$$\int_{0}^{1}\left<\frac{\partial v}{\partial s}, J\left(\frac{\partial
v}{\partial t} -X^{\alpha+dH_{t}}(v)\right)\right>dt=$$
$$=
\int_{0}^{1}\left<\frac{\partial v}{\partial s}, 
J(J_{t,s}\frac{\partial v}{\partial
s}+X^{\alpha+\chi(s)\sigma\eta+dH_{t}}(v)-X^{\alpha+dH_{t}}(v)\right>dt
=$$
$$=\int_{0}^{1}\left<\frac{\partial v}{\partial s}, 
JJ_{t,s}\frac{\partial v}{\partial s}\right>dt,$$
since $\eta|_{B}=0$ (so $X^{\chi(s)\sigma\eta}(v)=0$). It follows that :
$$\int_{0}^{1}\left<\frac{\partial v}{\partial s}, J(\frac{\partial
v}{\partial t} -X^{\alpha+dH_{t}}(v))\right>dt=$$
$$=-\int_{0}^{1}\left<J\frac{\partial v}{\partial s}, J_{s,t}\frac{\partial
v}{\partial s}\right>dt = -\int_{0}^{1}\left |\left|\frac{\partial
v}{\partial s}\right|\right|^{2}dt +\int_{0}^{1}\left<J\frac{\partial
v}{\partial s},(J_{s,t}-J)\frac{\partial
v}{\partial s}\right>dt\leq$$
$$\leq  -\int_{0}^{1}\left |\left|\frac{\partial
v}{\partial s}\right|\right|^{2}dt +\delta  \int_{0}^{1}\left |\left|\frac{\partial
v}{\partial s}\right|\right|^{2}dt= -(1-\delta)\left |\left|\frac{\partial
v}{\partial s}\right|\right|^{2}_{L^{2}}.$$
Choosing $\delta<2/3$ we get  : 
$$\int_{0}^{1}\left<\frac{\partial v}{\partial s}, J\left(\frac{\partial
v}{\partial t} -X^{\alpha+dH_{t}}(v)\right)\right>dt\, \leq\, -1/3
\left |\left|\frac{\partial
v}{\partial s}\right|\right|^{2}_{L^{2}}.\gol (4)$$

\noindent $2^{\circ}$ $Im(v(s,\cdot))\not\subset B$

Proceeding as above we obtain 
$$\int_{0}^{1}\left<\frac{\partial v}{\partial s}, J\left(\frac{\partial
v}{\partial t} -X^{\alpha+dH_{t}}(v)\right)\right>dt=$$
$$=
\int_{0}^{1}\left<\frac{\partial v}{\partial s}, 
J\left(J_{t,s}\frac{\partial v}{\partial
s}+X^{\alpha+\chi(s)\sigma\eta+dH_{t}}(v)-X^{\alpha+dH_{t}}(v)\right)\right>dt
=$$
$$=\int_{0}^{1}\left<\frac{\partial v}{\partial s}, 
JJ_{t,s}\frac{\partial v}{\partial s}\right>dt + \int_{0}^{1}\left<
\frac{\partial v}{\partial s}, JX^{\chi(s)\sigma\eta}(v)\right>dt=$$
$$=-\int_{0}^{1}\left<J\frac{\partial v}{\partial s}, J_{s,t}\frac{\partial
v}{\partial s}\right>dt -\int_{0}^{1}\omega_{M}\left(\frac{\partial v}{\partial s}, 
X^{\chi(s)\sigma\eta}(v)\right)dt=$$
$$= -\int_{0}^{1}\left |\left|\frac{\partial
v}{\partial s}\right|\right|^{2}dt +\int_{0}^{1}\left<J\frac{\partial
v}{\partial s},(J_{s,t}-J)\frac{\partial
v}{\partial s}\right>dt -\int_{0}^{1}
\chi(s)\sigma\eta\left(\frac{\partial v}{\partial s}\right)dt\leq $$
$$\leq  -(1-\delta)\left |\left|\frac{\partial
v}{\partial s}\right|\right|^{2}_{L^{2}} +\epsilon||\eta||\int_{0}^{1}\left|\left| 
\frac{\partial v}{\partial s}\right|\right|dt.$$
We apply the Cauchy Schwarz  inequality to the
 last integral and the fact that $\epsilon$
was chosen to satisfy $\epsilon ||\eta||<c/3$. We have therefore 
$$\int_{0}^{1}\left<\frac{\partial v}{\partial s}, J(\frac{\partial
v}{\partial t} -X^{\alpha+dH_{t}}(v))\right>dt\leq   -
(1-\delta)\left |\left|\frac{\partial
v}{\partial s}\right|\right|^{2}_{L^{2}}+\frac{c}{3}\left |\left|\frac{\partial
v}{\partial s}\right|\right|_{L^{2}}. \gol(5) $$

  We have not used the condition  $Im(v(s,\cdot))\not
\subset B$ yet. Note that it implies
using \ref{Palais} that 
$$\left|\left|\frac{\partial v}{\partial t}-X^{\alpha+dH_{t}}(v)\right|\right|\geq c.$$
 We infer :
$$\left|\left|\frac{\partial
v}{\partial s}\right|\right|_{L^{2}}=\left|\left|J_{s,t}\left(\frac{\partial
v}{\partial t}-X^{\alpha+\chi(s)\sigma\eta+dH_{t}}\right)\right|\right|_{L^{2}}\geq$$
$$\geq\left|\left|J\left(\frac{\partial
v}{\partial t}-X^{\alpha+\chi(s)\sigma\eta+dH_{t}}\right)\right|\right|_{L^{2}}-\, 
\left|\left|(J-J_{s,t})\left(\frac{\partial
v}{\partial t}-X^{\alpha+\chi(s)\sigma\eta+dH_{t}}\right)\right|\right|_{L^{2}}\geq$$
$$\geq
(1-\delta)\left|\left|\frac{\partial
v}{\partial t}-X^{\alpha+\chi(s)\sigma\eta+dH_{t}}\right|\right|_{L^{2}}\geq$$
$$\geq (1-\delta)\left(\left|\left|\frac{\partial
v}{\partial t}-X^{\alpha+dH_{t}}\right|\right|_{L^{2}}-\left|\left|
X^{\chi(s)\sigma\eta}\right|\right|_{L^{2}}\right)
\geq$$
$$\geq  (1-\delta)
\left(\left|\left|\frac{\partial
v}{\partial t}-X^{\alpha+dH_{t}}\right|\right|_{L^{2}}-\epsilon||\eta||\right)\geq$$
$$\geq (1-\delta)(c-c/3)=\frac{2}{3}(1-\delta)c, \gol (6)$$
using \ref{Palais}. 
 Let us come back to the inequality $(5)$. Using $(6)$ we find that for $\delta$ small
enough we have : $$ -
(1-\delta)\left |\left|\frac{\partial
v}{\partial s}\right|\right|^{2}_{L^{2}}+\frac{c}{3}\left |\left|\frac{\partial
v}{\partial s}\right|\right|_{L^{2}}\leq -1/3 \left |\left|\frac{\partial
v}{\partial s}\right|\right|^{2}_{L^{2}}.$$
Indeed, this is equivalent to 
$$\left(\frac{2}{3}-\delta\right)\left |\left|\frac{\partial
v}{\partial s}\right|\right|_{L^{2}}\geq c/3,$$
which, using $(6)$, is implied by
$$2(1-\delta) \left(\frac{2}{3}-\delta\right)\geq 1,$$
and this is true if we take for instance $\delta <1/10$. 

This means that in the case $2^{\circ}$ we also have 
$$\int_{0}^{1}\left<\frac{\partial v}{\partial s}, J\left(\frac{\partial
v}{\partial t} -X^{\alpha+dH_{t}}(v)\right)\right>dt\, \leq\, -1/3
\left |\left|\frac{\partial
v}{\partial s}\right|\right|^{2}_{L^{2}}\gol (7),$$
and the two cases imply that the above inequality is valid for any $s\in\real$. Now we
integrate this inequality with respect to $s$ and we get using $(3)$ that 
$$\hat{\cala}(g\bar{x}) -\hat{\cala}(\bar{y})\geq1/3
\int_{-\infty}^{+\infty}\left |\left|\frac{\partial
v}{\partial s}\right|\right|^{2}_{L^{2}[0,1]}ds=\frac{1}{3 }E(v).$$
This implies the inequality $(2)$ and finishes the proof of \ref{crucial}. \hfill
$\diamond$

  Therefore the space $\calm^{0}_{g,s}(x,y)$ is finite. Moreover, if  $m_{g}(x,y)$ is
its parity the relation $(2)$ in the proof of \ref{crucial} shows that
the sum $$\sum_{g\in\piu(M)}m_{g}(x,y)g$$ belongs to $\Lambda_{u}$. This enables us to
define the morphism $$\Gamma_{\bullet} :    C_{\bullet}(L, \phi_{t}^{\alpha+dH_{t}},
J_{t})\ri C_{\bullet}(L, \psi_{t}^{\alpha+\sigma\eta+dH_{t}}, J'_{t}).$$
The fact that $\Gamma_{\bullet}$ commutes with the differentials is proved by studying
the compactness of the $1$-dimensional components $\calm^{1}_{g,s}(x,y)$. The proof is
analogous to \ref{bord}. To show that $\Gamma_{\bullet}$ induces an isomorphism at the
level of homology, one should proceed as usual : Take a  homotopy between $(\alpha
+\sigma\eta+dH_{t},J'_{t})$  and $(\alpha+dH_{t},J_{t})$. This yields a morphism
$$\Gamma'_{\bullet} :
C_{\bullet}(L, \psi_{t}^{\alpha+\sigma\eta+dH_{t}}, J'_{t})\ri
C_{\bullet}(L, \phi_{t}^{\alpha+dH_{t}}, J_{t}).$$
Then, using two-parameter homotopies as in \cite{FS} one can show  that
$\Gamma_{\bullet}\Gamma'_{\bullet}$ is homotopic to the identity and also that  
$\Gamma'_{\bullet}\Gamma_{\bullet}$ is homotopic to the identity. This shows that the
homologies of the complexes are isomorphic, and the proof of \ref{epsilon2}.a is
finished.\\

\noindent \underline{Proof of \ref{epsilon2}.b}

Replace overall in the proof of \ref{epsilon2}.a the $1$-form $\alpha\in u$ with an
exact $1$-form $df$ where $f:M\ri \real$. We get then  in a similar manner 
  a $\Lambda_{u}$-morphism
$$\Gamma^{0}_{\bullet}:C_{\bullet}(L, \phi_{t}^{df+dH_{t}}, J_{t})\ri
C_{\bullet}(L, \phi_{t}^{df+\sigma\eta + dH_{t}}, J'_{t}),$$
which is an isomorphism in homology. The first complex above is actually a
$\Lambda$-complex with coefficients extended to $\Lambda_{u}$. We can replace it by
$\Lambda_{u}\otimes_{\Lambda}C_{\bullet}(L, \phi_{t}^{df+dH_{t}}, J_{t})$. Now the
Hamiltonian invariance implies that the homotopy class of the $\Lambda$-complex
$C_{\bullet}(L, \phi_{t}^{df+dH_{t}}, J_{t})$ does not depend on the choice of a
regular pair $(df+dH_{t}, J_{t})$. On the other hand, the classical Floer argument
\cite{Flo3} provides  a regular pair $(H^{0}_{t}, J^{0}_{t})$ such that
the complex $C_{\bullet}(L, \phi_{t}^{dH^{0}_{t}}, J^{0}_{t})$ coincides 
with the Morse
complex defined by a pair $(g, \xi)$ --~ where $g:L\ri\real$ is a Morse function 
and $\xi$ is a
generic pseudogradient on $L$ -- in the following sense : the complex
$C_{\bullet}(\phi_{t}^{dH^{0}_{t}}, J^{0}_{t})$ is spanned by the
intersection points $L\cap (L+dg)$ in $\tastl\subset \tast$ and the map 
$$v(s,t)\mapsto v(s,0)$$ is a one-to-one correspondance between the holomorphic 
strips 
which define the Floer differential and the trajectories of $\xi$ which define the
Morse differential. 

It is then easy to see that the above Morse complex is identical to
$C_{\bullet}~(~\bar{L}~\ri~ L~, ~\xi~)$, where $\bar{L}\ri L$ is the pull-back of
$\tastilde\ri\tast$ : this is the Morse complex defined using $(g,\xi)$, by lifting
the trajectories of $\xi$ to $\bar{L}$, as in  Section 2.2 (see Remark
\ref{novikovgeneral}).
Moreover, the Novikov ring which defines the Novikov homology associated to the class
$p^{\ast}u$ and to the covering $\bar{L}\ri L$ is the same as $\Lambda_{u}$. Finally,
if we denote by $\sim$ the relation of 
homotopy equivalence between $\Lambda_{u}$-complexes we
get : 
$$\Lambda_{u}\otimes_{\Lambda}C_{\bullet}(L, \phi_{t}^{df+dH_{t}}, J_{t})\sim
\Lambda_{u}\otimes_{\Lambda}C_{\bullet}(L, \phi_{t}^{dH^{0}_{t}}, J^{0}_{t})\sim
\Lambda_{u}\otimes_{\Lambda}C_{\bullet}(\bar{L}\ri L, g, \xi).$$
  The latter complex defines the Novikov homology 
$H(\bar{L}\ri  L, p^{\ast}u)$ corresponding to
the covering $\bar{L}\ri L$. We therefore have the isomorphism :
$$\bar{FH}(L,\sigma u)\, \approx\, H(\bar{L}\ri L,p^{\ast}u),$$
 and the proof of \ref{epsilon2}.b. is complete.
\hfill $\diamond$

\vspace{.2in}

\noindent\underline{End of the proof of \ref{epsilon}}

We argue as in Subsection 3.7.  The holomorphic strips of the moduli space
$\calm^{0}_{g, s}(x,y)$ define paths on $L$ from $x$ to $y$ by the formula : 
$$v\, \mapsto \, v(s,0).$$
For fixed lifts $\tilde{x}, \, \tilde{y}\, \in\, \widetilde{L}$, such a path lifts
to $\widetilde{L}$ to a path joining $h\tilde{x}$ and $\tilde{y}$ for some $h\in
\piu(L)$. We define thus the spaces $ \widetilde{\calm}^{0}_{h, s}(x,y)$. We
obviously have :
$$\calm^{0}_{g, s}(x,y)\, =\, \bigcup_{h\in \piu(L), \, p(h)=g}
\widetilde{\calm}^{0}_{h, s}(x,y).$$
In particular, any set $\widetilde{\calm}^{0}_{h, s}(x,y)$ has a finite number of
elements  (by \ref{crucial}) and the formula 
$$\widetilde{\Gamma}(x)\, =\, \sum_{h\in \piu(L), y}m_{h}(x,y)hy$$
 -- where $m_{h}(x,y)$ is the number modulo $2$ of elements in
$\widetilde{\calm}^{0}_{h, s}(x,y)$ -- defines a morphism of
$\Lambda_{p^{\ast}u}$-complexes
$$\widetilde{\Gamma}_{\bullet}:   
\widetilde{C}_{\bullet}(L, \phi_{t}^{\alpha+dH_{t}},
J_{t})\ri \widetilde{C}_{\bullet}(L, \psi_{t}^{\alpha+\sigma\eta+dH_{t}}, J'_{t})$$
which is a lift of the morphism $\Gamma_{\bullet}$ defined in the proof of
\ref{epsilon2}. Using the same argument we show that it induces an isomorphism in
homology, proving thus \ref{epsilon}.a. Then, we proceed analogously in order to
define a lift of the morphism $\Gamma_{\bullet}^{0}$ of the proof of
\ref{epsilon2}.b. For the particular choice $(H_{t}^{0}, J_{t}^{0})$ of
\ref{epsilon2}.b, the lifted complex $\widetilde{C}_{\bullet}$ is the Morse complex
defined on the universal cover $\widetilde{L}$. Therefore, the same argument as in
the proof of \ref{epsilon2}.b yields
an isomomorphism 

$$HF(L,\sigma u) \, \approx\, H(L,p^{\ast}u),$$
completing the proof of \ref{epsilon}. 

\hfill $\diamond$

As we explained, this immediately implies \ref{iso}. \hfill $\diamond$

\subsection{Proof of the main results}

\noindent\underline{Proof of \ref{novi}}\\

If $M$ is the total space of a fibration over the circle and $L\subset\tast $ is
exact, we consider a non-vanishing closed $1$-form $\alpha$ on $M$ and we define the
symplectic isotopy $\Psi_{t}(p,q)=(p+t\alpha_{q}, q)$.  Obviously $\Psi_{T}(L)\cap
L=\emptyset$ for $T$ sufficiently large. Without restricting the generality we may
suppose that $T=1$. The Floer complex defined in the previous subsection is empty so
$FH(L,u)=0$, where $u$ is the cohomology class of $\alpha$. Using \ref{iso} this
implies \ref{novi}, i.e. :
 $$H_{\ast}(L, p^{\ast}u)=0. \gol(1)$$ Note that since $p:\piu(L)\ri \piu(M)$
 is 
an epimorphism the class 
$p^{\ast}u$ is not zero. 

\hfill $\diamond$

\noindent\underline{Proof of \ref{main}}\\

 a) Let $<g_{1}, g_{2}, \ldots, g_{p}|r_{1},\ldots, r_{q}>$ be a presentation of
$\piu(L)$. If $p-q\geq 2$ then by \ref{nonzero}.a we vave that 
 $ H_{1}(L, p^{\ast}u)\neq 0$, contradicting $(1)$. 
We infer that  $p-q\leq 1$, as claimed. \\
\\
b) If $\piu(L)=G_{1}\ast G_{2}$ for some non trivial groups $G_{i}$, then we have
again  $ H_{1}(L, p^{\ast}u)\neq 0$, by \ref{nonzero}.b, and the proof is finished.

\hfill $\diamond$

\vspace{.2in} 

\noindent\underline{Proof of \ref{coro}}\\
\\
a) Suppose that there exists an exact Lagrangian embedding $$L\times P\hookrightarrow
T^{\ast}(Q\times{\bf S}^{1}),$$ where $\chi(L)\neq 0$ and $\piu(P)$ is finite. 
As above, we obtain $H_{\ast}(L\times P, p^{\ast}u)= 0$, where $u$ is the
class of the nonvanishing closed $1$-form $d\theta$ on $Q\times{\bf S}^{1}$. Since
$\piu(P)$ is finite, we obtain $ p^{\ast}u \in H^{1}(L, \real)\subset  H^{1}(L\times
P, \real)$. 
We  apply  \ref{Kunneth} and  we obtain $H_{\ast}(L, p^{\ast}u)= 0$. But this
contradicts \ref{nonze}.\\ 
\\
b)    Suppose that there exists a Lagrangian embedding $$(L_{1}\#L_{2})
\times P\hookrightarrow
T^{\ast}(Q\times{\bf S}^{1}).$$
We show that either $L_{1}$ or $L_{2}$ is a simply connected ${\bf Z}/2$-homology
sphere. As $n \geq
4$, the fundamental group of $L=L_{1}\#L_{2}$ is  the free product $\piu(L_{1})\ast
\piu(L_{2})$. We get from \ref{main} that one of the $L_{i}$'s is simply connected.
Suppose $\piu(L_{1})=1$. 

By \ref{novi} we know that $H_{\ast}(L, p^{\ast}u)=0$, where $u$ is the
class  $d\theta$ on $M=Q\times{\bf S}^{1}$ and $p$ is the projection. 
 We show that this Novikov homology  cannot vanish unless $L_{1}$ is
a simply connected ${\bf Z}/2$-homology sphere.

 Denote $D_{i}\subset L_{i}$ two embedded open $n$-disks
and write $L$ as 
$$(L_{1}\setminus D_{1})\, \bigcup_{{\bfs}^{n-1}}\, (\bfs^{n-1}\times [0,1])\,  
\bigcup_{{ \bfs}^{n-1}}\, (L_{2}\setminus D_{2}).$$
Choose a $CW$-structure on $L$ which fits to this decomposition and which is
the standard product structure
on $\bfs^{n-1}\times [0,1]$ (we take the decomposition with one $0$-cell and one
$(n-1)$-cell on $\bfs^{n-1}$). Denote by $D_{\bullet}$ the $\bf Z$$/2$-free 
subcomplex spanned by the cells of 
$L_{1}\setminus D_{1}$. Consider the Novikov complex $C_{\bullet}(L,p^{\ast}u)$ which
is by definition the tensor product $\Lambda_{p^{\ast}u}\otimes_{{\bf
Z}/2[\piu(L)]}C_{\bullet}(\widetilde{L})$. Since $L_{1}$ is simply connected,
 $\widetilde{L}$ is the connected sum of $\widetilde{L}_{2}$ with $\piu(L)$ copies
of $L_{1}$. Using this fact and our  choice of the $CW$-structure, we get
 that   the complex
$$\Lambda_{p^{\ast}u}\otimes_{{\bf Z}/2}D_{1\leq\bullet\leq n-1}$$
is a direct summand in $C_{\bullet} (L,p^{\ast}u)$. In particular, since the Novikov
complex is acyclic (by \ref{novi}), the homology of $D_{\bullet}$, which is
 the 
homology $H_{i}~(L_{1}~\setminus~D_{1}~;{\bf Z}/2)$, vanishes in degrees
 $2\leq i\leq n-2$.
Then, $L_{1}\setminus D_{1}$ is simply connected so the homology 
$H_{1}(L_{1}\setminus D_{1})$ vanishes, too. Using
Poincar\'e duality, we  find that
$$H_{n-1}(L_{1}\setminus D_{1})\approx H^{1}(L_{1}\setminus D_{1},\partial D_{1})$$
and the latter vanishes since $L{1}\setminus D_{1}$ is simply connected. We also have : 
$$H_{n}(L_{1}\setminus D_{1})\approx H^{0}(L_{1}\setminus D_{1},\partial D_{1})=0.$$
Therefore,  the groups $H_{n-1}$ and $H_{n}$ are 
also zero,
so $L_{1}\setminus D_{1}$ has the ${\bf Z}/2$ homology of the $n$-disk.  

Using Mayer-Vietoris we find that $L_{1}$ is a ${\bf Z}/2$-homology sphere. \\
\\
c) Denote by  $H\subset H^{1}(T^{m}\times Q; {\bf Z})$ the subgroup $H^{1}(T^{m}\times
\{pt\};{\bf Z})$. The group $H$ is isomorphic to ${\bf Z}^{m}$ 
and any cohomology class $u\in
H\setminus\{0\}$ obviously contains a closed non-vanishing $1$-form on $T^{m}\times
Q$. It follows by  \ref{novi} that $$H_{\ast}(L\times T^{l}, p^{\ast}u)=0$$
for every $u\in H\setminus\{0\}$. 

Recall that we may suppose (using \ref{epi}) that $p$
is an epimorphism, which implies that $p^{\ast}$ is a monomorphism.  
Consider  $H^{1}(L)\subset
H^{1}(L\times T^{l})$. We show that 
$$H^{1}(L)\cap p^{\ast}(H)\, \neq\, \{0\}.$$
If not, then the composition $$ {\bf Z}^{m}\approx p^{\ast}H \hookrightarrow 
H^{1}(L\times T^{l}) \stackrel{pr}{\ri} H^{1}(T^{l})\approx{\bf Z}^{l}$$
is a monomorphism, which is impossible, since $l<m$. We infer that there exists a
nonvanishing class
$p^{\ast}u=v\in
H^{1}(L)\subset
H^{1}(L\times T^{l})$ such that $$H_{\ast}(L\times T^{l},v)=0.$$ By applying
\ref{Kunneth} we obtain then $H_{\ast}(L,v)=0$. Finally, we apply \ref{nonzero} to get
the desired conclusions on $\piu(L)$ and finish the proof.
\hfill $\diamond$

\vspace{.2in}

\noindent {\bf Acknowledgements}  I thank Michele Audin for many valuable discussions
on Floer homology. I thank Agnes Gadbled and Alexandru Oancea, who carefully read
the manuscript and helped me to improve it by their remarks and corrections. Finally, I thank the
anonimous referee for his useful remarks and suggestions.

\end{document}